\documentclass[oneside,english,reqno]{amsart}
\usepackage[T1]{fontenc}
\usepackage[latin9]{inputenc}
\usepackage{xcolor}
\usepackage{pdfcolmk}
\usepackage{babel}
\usepackage{float}
\usepackage{mathrsfs}
\usepackage{amsthm}
\usepackage{amstext}
\usepackage{amssymb}
\usepackage{graphicx}
\usepackage{esint}
\usepackage[all]{xy}
\PassOptionsToPackage{normalem}{ulem}
\usepackage{ulem}
\usepackage[unicode=true,pdfusetitle,
 bookmarks=true,bookmarksnumbered=false,bookmarksopen=false,
 breaklinks=false,pdfborder={0 0 1},backref=false,colorlinks=true]
 {hyperref}
\hypersetup{
 colorlinks=true,citecolor=blue,linkcolor=blue,linktocpage=true}

\makeatletter

\newcommand{\lyxmathsym}[1]{\ifmmode\begingroup\def\b@ld{bold}
  \text{\ifx\math@version\b@ld\bfseries\fi#1}\endgroup\else#1\fi}

\providecommand{\tabularnewline}{\\}
\newcommand{\lyxdot}{.}

\providecolor{lyxadded}{rgb}{0,0,1}
\providecolor{lyxdeleted}{rgb}{1,0,0}
\numberwithin{equation}{section}
\numberwithin{figure}{section}
\theoremstyle{plain}
\newtheorem{thm}{\protect\theoremname}[section]
  \theoremstyle{plain}
  \newtheorem{lem}[thm]{\protect\lemmaname}
  \theoremstyle{definition}
  \newtheorem{defn}[thm]{\protect\definitionname}
  \theoremstyle{plain}
  \newtheorem{prop}[thm]{\protect\propositionname}
  \theoremstyle{plain}
  \newtheorem{cor}[thm]{\protect\corollaryname}
  \theoremstyle{remark}
  \newtheorem{rem}[thm]{\protect\remarkname}
  \theoremstyle{definition}
  \newtheorem{example}[thm]{\protect\examplename}

\providecommand{\MR}[1]{}

\makeatother

  \providecommand{\corollaryname}{Corollary}
  \providecommand{\definitionname}{Definition}
  \providecommand{\examplename}{Example}
  \providecommand{\lemmaname}{Lemma}
  \providecommand{\propositionname}{Proposition}
  \providecommand{\remarkname}{Remark}
\providecommand{\theoremname}{Theorem}

\begin{document}

\title[Momentum Operators in Two Intervals]{Momentum Operators in Two Intervals: \\ Spectra and Phase Transition}

\author{Palle E. T. Jorgensen, Steen Pedersen, and Feng Tian}

\address{(Palle E.T. Jorgensen) Department of Mathematics, The University
of Iowa, Iowa City, IA 52242-1419, U.S.A. }

\email{jorgen@math.uiowa.edu }

\urladdr{http://www.math.uiowa.edu/\textasciitilde{}jorgen/}

\address{(Steen Pedersen) Department of Mathematics, Wright State University,
Dayton, OH 45435, U.S.A. }

\email{steen@math.wright.edu }

\urladdr{http://www.wright.edu/\textasciitilde{}steen.pedersen/}

\address{(Feng Tian) Department of Mathematics, Wright State University, Dayton,
OH 45435, U.S.A. }

\email{feng.tian@wright.edu }

\urladdr{http://www.wright.edu/\textasciitilde{}feng.tian/}
\begin{abstract}
We study the momentum operator defined on the disjoint union of two
intervals. Even in one dimension, the question of two non-empty open
and non-overlapping intervals has not been worked out in a way that
extends the cases of a single interval and gives a list of the selfadjoint
extensions. Starting with zero boundary conditions at the four endpoints,
we characterize the selfadjoint extensions and undertake a systematic
and complete study of the spectral theory of the selfadjoint extensions.
In an application of our extension theory to harmonic analysis, we
offer a new family of spectral pairs. Compared to earlier studies,
it yields a more direct link between spectrum and geometry.
\end{abstract}


\maketitle
\tableofcontents{}

\section{\label{sec:Introduction}Introduction}

An early paper of John von Neumann \cite{vNeu32} deals with unbounded
linear operators in Hilbert space, and their spectral theory. By now
the theory is part of functional analysis, and the problems entail
the study of \textquotedblleft{}deficiency indices;\textquotedblright{}
now important tools in mathematical physics; see e.g., \cite{Kr49}.

Indeed, von Neumann\textquoteright{}s initial motivation came directly
from quantum mechanics: measurements of observables, prepared in states.
These indices have now become the accepted framework in spectral theory,
as well as in quantum mechanics.

In quantum theory, states are realized by vectors (of norm one) in
Hilbert space $\mathscr{H}$, and the quantum observables correspond
to selfadjoint operators in $\mathscr{H}$ (the state space.)

Our system of operators in two intervals has a number of applications
not contained in other studies. For example, when modeling quantum
states (Section \ref{sec:U(T)}), we show how a two-component quantum
system offer a binary model for the dynamics of wave functions. Indeed
the parameters for our family of selfadjoint operators allow a direct
computation of transition probabilities. 

Here we study an instance of the index/spectral theory question, and
in computations we use the spectral theorem as it is applied to selfadjoint
operators. However the deficiency indices are precisely the obstruction
to selfadjointness for unbounded Hermitian operators.

This is illustrated with examples, naturally occurring, and quite
simple, from quantum theory: There one studies wave-functions with
confinement to a finite region (particles in a box), e.g., confinement
in a finite interval. In these cases, the physical operators are not
selfadjoint, but only Hermitian.

But the unbounded linear operators will still have dense domain in
the physical Hilbert space (of states). Therefore, an assignment of
boundary conditions will amount to a subtle extension of the initial
dense operator domain.

In several variables, one is then looking for commuting selfadjoint
extensions; see \cite{Fu74,Jo82,PW01,Pe04,JP00,JP99,JP98,DJ07} for
a study of this problem, in several applications.

For these problems, and for applications in the theory of linear partial
differential operators, one is led to a more axiomatic consideration
of boundary conditions. Indeed, J. von Neumann, in his paper \cite{vNeu32},
identified such conditions in the Hilbert space theoretic framework.
For a single Hermitian operator, they take the form of deficiency
spaces, and deficiency indices, see details below. In one dimension,
and for a finite interval, say $I$, a popular model example has its
deficiency indices come out $(1,1).$ Hence, by von Neumann\textquoteright{}s
theorem, this yields a corresponding one-parameter family of selfadjoint
extensions. In the case of one interval, they can readily be computed,
and it is possible to read off the spectrum of each of the selfadjoint
extensions. Each one is the generator of a unitary one-parameter group
$U(t)$ acting in $\mathscr{H}$. And each of these unitary groups
of operators in turn translates the wave-functions in the interior
of the interval, and then assigns a phase-factor at the endpoints
of $I$. The corresponding dynamics of the wave-functions turns into
a fixed phase-shift, i.e., a matching with a phase-factor of incoming
and outgoing waves.

Moreover, for each selfadjoint extension, the corresponding spectrum
is a Fourier spectrum for the Hilbert space $\mathscr{H}=L^{2}(I)$
where $I$ is the given interval. By \textquotedblleft{}Fourier spectrum\textquotedblright{}
we mean that the spectrum (of the selfadjoint extension) is discrete
of uniform multiplicity one, and that the eigenfunctions are complex
exponentials. We say that the pair, set and spectrum, is a spectral
pair. See \cite{JP98} for the theory of spectral pairs.

Indeed, for more general models with wave equations, and in higher
dimensions, there is a scattering theory model which accounts for
scattering of waves at obstacles, see for example \cite{LP68}. Yet,
still staying with the translation model, but going from one dimension
to $\mathbb{R}^{d},$ Fuglede \cite{Fu74} asked for an analogue of
the translation model in one dimension, but now for bounded open domains
in $\mathbb{R}^{d}.$ 

Nonetheless, even for one dimension, the question of two non-empty
open and non-overlapping intervals has not been worked out in a way
that extends the cases of a single interval and gives a list of the
selfadjoint extensions.

\subsection{Two Intervals}

This paper is devoted to the two-interval case. Now, starting with
zero boundary conditions at the four endpoints, we show that the deficiency
indices are then $(2,2),$ see \cite{Kr49}. Moreover, we offer a
characterization and parameterization of all the selfadjoint extensions.
For each one, we compute the spectrum, and we write down the corresponding
unitary one-parameter groups.

We show that, by contrast to the one-interval case, that some of the
selfadjoint extension may have points in the spectrum of multiplicity
two. 

Our analysis lies at the crossroads of operator theory, spectral geometry,
and harmonic analysis. For the two-interval case, we find the selfadjoint
extension operators $P_{B},$ and their spectra. Each operator $P_{B},$
indexed by the group $U(2)$ of unitary $2$ by $2$ matrices, has
pure eigenvalue spectrum $\{\lambda_{n}(B)\};$ the multiplicity can
be $1$ or $2.$ In fact, we have the eigenvalues represented as functions
of both $B$ and of the choice of endpoints in $I_{2}.$

\textcolor{black}{We study all possible configurations of the two
intervals. It turn out that only the ratio $r$ of their lengths is
important. }The nature of this number has subtle spectral theoretic
implications; e.g., the distinction between the cases when $r$ is
rational vs irrational is important. When it is irrational, the points
in the spectrum of $P_{B}$ form a generator for pseudo-random numbers
\cite{BAAS11}; see Section \ref{sec:Spectral-Theory} for details.
For the general case, we display the points $\{\lambda_{n}\}$ in
the spectrum of $P_{B}$ as functions of 5 variables, 4 for $B$ in
$U(2),$ and one for $r.$

In Section \ref{sec:Spectral-Pairs}, we determine completely the
cases when we the two intervals, together with spectrum of $P_{B},$
yield spectral pairs (SP). It turns out to be a small subset of the
total parameter space. In the SP-case, we show that the set $\Lambda$
of points $\left\{ \lambda_{n}\right\} $ fall in two simple cases:
Modulo the angular part of $B,$ $\Lambda$ is either is the union
of one or two lattice co-sets. For simplicity suppose the intervals
are $[0,1]\cup[\alpha,\beta]$ where $1<\alpha<\beta.$ In Theorems
\ref{thm:-spectral-w=00003D0} and \ref{thm:spectral-0<w<1}, we give
a complete analysis of the SP cases. Modulo angular variables, this
is merely a \textquotedblleft{}small\textquotedblright{} discrete
subfamily: Firstly, SP forces either $\beta/(1+\beta-\alpha)$ to
be an integer or $\beta-\alpha=1;$ and secondly, we get the eigenvalue
lists as follows: Let $w$ be the radial variable in $B.$ If $\beta/(1+\beta-\alpha)$
is an integer then the points $\lambda_{n}$ form a lattice, and if
$\beta-\alpha=1$ the points $\lambda_{n}$ (as a function of $w$)
contains the branches of the arc-cos function. In summary, the SP
cases form a singular family, representing a broken symmetry, or phase
transition; a discrete family within a much bigger continuum. For
comparison, note that in \cite{Rob71}, dealing with second order
operators in a single interval, the author finds an analogous spectral
analysis of phase transitions.

\subsection{Unbounded Operators}

We recall the following fundamental result of von Neumann on extensions
of Hermitian operators.
\begin{lem}[see e.g. \cite{DS88b}]
\label{lem:vN def-space}Let $L$ be a \textcolor{black}{closed}
Hermitian operator with dense domain $\mathscr{D}_{0}$ in a Hilbert
space. Set 
\begin{eqnarray}
\mathscr{D}_{\pm} & = & \{\psi_{\pm}\in dom(L^{*})\left.\right|L^{*}\psi_{\pm}=\pm i\psi_{\pm}\}\nonumber \\
\mathscr{C}(L) & = & \{U:\mathscr{D}_{+}\rightarrow\mathscr{D}_{-}\left.\right|U^{*}U=P_{\mathscr{D}_{+}},UU^{*}=P_{\mathscr{D}_{-}}\}\label{eq:vN1}
\end{eqnarray}
where $P_{\mathscr{D}_{\pm}}$ denote the respective projections.
Set 
\[
\mathscr{E}(L)=\{S\left.\right|L\subseteq S,S^{*}=S\}.
\]
Then there is a bijective correspondence between $\mathscr{C}(L)$
and $\mathscr{E}(L)$, given as follows: 

If $U\in\mathscr{C}(L)$, and let $L_{U}$ be the restriction of $L^{*}$
to 
\begin{equation}
\{\varphi_{0}+f_{+}+Uf_{+}\left.\right|\varphi_{0}\in\mathscr{D}_{0},f_{+}\in\mathscr{D}_{+}\}.\label{eq:vN2}
\end{equation}
Then $L_{U}\in\mathscr{E}(L)$, and conversely every $S\in\mathscr{E}(L)$
has the form $L_{U}$ for some $U\in\mathscr{C}(L)$. With $S\in\mathscr{E}(L)$,
take 
\begin{equation}
U:=(S-iI)(S+iI)^{-1}\left.\right|_{\mathscr{D}_{+}}\label{eq:vN3}
\end{equation}
and note that 
\begin{enumerate}
\item $U\in\mathscr{C}(L)$, and
\item $S=L_{U}$.
\end{enumerate}

Vectors $f$ in $dom(L^{*})$ admit a unique decomposition $f=\varphi_{0}+f_{+}+f_{-}$
where $\varphi_{0}\in dom(L)$, and $f_{\pm}\in\mathscr{D}_{\pm}$.
For the boundary-form $\mathbf{B}(\cdot,\cdot)$, we have
\begin{eqnarray*}
i\mathbf{B}(f,f) & = & \left\langle L^{*}f,f\right\rangle -\left\langle f,L^{*}f\right\rangle \\
 & = & \left\Vert f_{+}\right\Vert ^{2}-\left\Vert f_{-}\right\Vert ^{2}.
\end{eqnarray*}

\end{lem}

\subsection{Summary}

We undertake a systematic study of interconnections between geometry
and spectrum for a family of selfadjoint operator extensions indexed
by two thing: by (i) the prescribed configuration of the two intervals,
and by (ii) the von Neumann parameters. This turns out to be subtle,
and we show in detail how variations in both (i) and (ii) translate
into explicit spectral properties for the extension operators. Indeed,
for each choice in (i), i.e., relative length of the two intervals,
we have a Hermitian operator with deficiency indices $(2,2).$ Our
main theme is spectral theory of the corresponding family of $(2,2)$-selfadjoint
extension operators. Accidentally, this more general problem has implications
for a question from \cite{Laba01} on spectral pairs. But our study
of the spectral theory dictated by (i) and (ii) goes far beyond its
implications for spectral pairs.

Another advantage of our approach via a definite $6$-parameter family
of selfadjoint operators is that the verification for completeness
of the associated function-systems becomes easy. We first prove that
all these operators have pure point-spectrum, see Theorem \ref{thm:Magic-equation}.
In each verification, one may then take advantage of the spectral
theorem applied to the particular selfadjoint operators under study.
So this applies whether or not one has a spectral pair. Otherwise
the completeness property is typically subtle in spectral geometry
and in physics applications; see e.g., \cite{DS88b,Gil72,KM06,Laba01,La02,LP68,Rob71}. 

One of our applications will be to the characterization of spectral
pairs in the two-interval case(s). While this question was addressed
in the paper \cite{Laba01}, it was attacked there with the use of
geometric tools; specifically by listing the ways two intervals can
tile the real line under translations. We obtain the same results
regarding spectral pairs, but now with the use of operator/spectral
theory. Even though most of our computations are for interval examples,
we believe they are of general interest, and that they offer intriguing
extensions to more structured models.

\subsection{Prior Literature}

There are related investigations in the literature on spectrum and
deficiency indices. For the case of indices $(1,1)$, see for example
\cite{ST10}. For a study of odd-order operators, see \cite{BH08}.
Operators of even order in a single interval are studied in \cite{Oro05}.
The paper \cite{BV05} studies matching interface conditions in connection
with deficiency indices $(m,m)$. Dirac operators are studied in \cite{Sak97}.
For the theory of selfadjoint extensions operators, and their spectra,
see \cite{Smu74,Gil72}, for the theory; and \cite{Naz08,VGT08,Vas07,Sad06,Mik04,Min04}
for recent papers with applications. For applications to other problems
in physics, see e.g., \cite{PoRa76,Ba49}.

\section{\label{sec:Momentum-Operators}Momentum Operators}

By momentum operator $P$ we mean the generator for the group of translations
in $L^{2}(-\infty,\infty)$. There are several reasons for taking
a closer look at restrictions of the operator $P.$ In our analysis,
we study spectral theory determined by a choice of two intervals (details
below.) Our motivation derives from quantum theory, and from the study
of spectral pairs in geometric analysis; see e.g., \cite{DJ07}, \cite{Fu74},
\cite{JP99}, \cite{Laba01}, and \cite{PW01}. In our model, we examine
how the spectral theory depends on both variations in the choice of
the two intervals, as well as on variations in the von Neumann parameters. 

Granted that in many applications, one is faced with vastly more complicated
data and operators; nonetheless, it is often the case that the more
subtle situations will be unitarily equivalent to a suitable model
involving $P$. This is reflected for example in the conclusion of
the Stone-von Neumann uniqueness theorem: The Weyl relations for quantum
systems with a finite number of degree of freedom are unitarily equivalent
to the standard model with momentum and position operators $P$ and
$Q$.

\subsection{Applications\label{sub:Applications}}

In Lax-Phillips scattering \cite{LP68} theory for the acoustic wave
equation, one obtains the solution to a particular wave equation represented
by a unitary one-parameter group acting on Hilbert space. Wave equations
preserve energy; hence the Hilbert space. 

With scattering against a finite obstacle in some number $d$ of dimensions,
one then gets two closed subspaces corresponding to incoming and outgoing
waves, respectively, here referring to the whole space, in this case
$\mathbb{R}^{d}$. The Lax-Phillips scattering operator is derived
from the two. Again, in the two subspaces, the respective unitary
one-parameter groups are unitary equivalent to a $P$-model. By taking
into account uniform multiplicity, one arrives at the operator $P$. 

For other uses of $P$ and its restrictions, see \cite{Jo81} and
\cite{Rob71}. The application in \cite{Rob71} is to quantum statistical
mechanics. 

In geometric analysis, one studies second order Hermitian PDOs on
a finite region in $\mathbb{R}^{d}$, and selfadjoint boundary condition,
for example Dirichlet conditions, or Neumann conditions, see also
\cite{LP68}, \cite{DS88b} and \cite{Rob71}. 

In quantum systems, the initial interest is in the Schr\"{o}dinger
equation for the dynamics of quantum states. And once again an essential
part of the problem is turned into the study of a unitary one-parameter
group in Hilbert space governing the solution to the Schrödinger equation:
We are faced with determining the spectrum of the selfadjoint generator.
Even if there may not be a direct unitary equivalence, there will
often be a part of the dynamics, e.g., bound-states, localization,
and quantum tunneling, where the two-interval model offers insight.
In \cite{Rob71} for example, the operator is second order, and it
is used in the study of confinement quantum states (wave functions)
in a single interval. There is a connection between this, and the
details below. Below we have a first-order operator, and two intervals.
But in both instances, the models are represented through the study
of deficiency indices $(2,2)$. In both cases, one is able to determine
an eigenvalue spectrum (i.e., bound states) from a system of curves
with intersection, but the systems vary from one application to the
other. For details, in \cite{Rob71}, see Figures 1 through 8. The
analogous system of curve intersections of our present two-interval
models are given in Examples \ref{Example: dense orbits} and \ref{examples:ratinal vs irrational}
below.

In studies of quantum scales, one is faced with quantum tunneling,
see e. g., \cite{Rob71}, a process that cannot be directly perceived
at a macroscopic scale. Quantum particles (or rather waves) travel
between potential barriers. With a small probability, wave functions
will crossing barriers. See Corollary \ref{cor:U(t)-1} for a discussion.
The reason is duality for quantum states; i.e., states having simultaneously
properties of waves and particles. This leads to computation of probability
density for particle's position, thus describing the probability that
the particle is at any given place. In the limit of large barriers,
the probability of tunneling decreases for taller and wider barriers.
The simplest and best understood tunneling-barriers are rectangular
barriers, offering approximate solutions.

\subsection{The boundary form, spectrum, and the group $U(2)$}

Since the problem is essentially invariant under affine transformations
we may assume the two intervals are $I_{1}:=[0,1]$ and $I_{2}:=[\alpha,\beta]$.
$L^{2}(I_{1}\cup I_{2})$ is a Hilbert space with respect to the inner
product
\begin{equation}
\langle f\mid g\rangle:=\int_{I_{1}}f\overline{g}+\int_{I_{2}}f\overline{g}.\label{eq:InnerProduct}
\end{equation}
 The \emph{maximal momentum operator} is 
\begin{equation}
P:=\frac{1}{i2\pi}\frac{d}{dt}\label{eq:MomentumOperator}
\end{equation}
with domain $\mathscr{D}(P)$ equal to the set of absolutely continuous
functions on $I_{1}\cup I_{2}.$ 

The \emph{boundary form} associated with $P$ is the form 
\begin{equation}
\mathbf{B}(f,g):=\langle Pf\mid g\rangle-\langle f\mid Pg\rangle\label{eq:BoundaryForm1}
\end{equation}
on $\mathscr{D}(P).$ Clearly, 
\begin{equation}
\mathbf{B}(f,g)=f(1)\overline{g(1)}-f(0)\overline{g(0)}+f(\beta)\overline{g(\beta)}-f(\alpha)\overline{g(\alpha)}.\label{eq:BoundaryForm2}
\end{equation}
For $f\in\mathscr{D}\left(P\right),$ let $\rho_{1}(f):=\left(f(1),f(\beta)\right)$
and $\rho_{2}(f):=\left(f(0),f(\alpha)\right).$ Then 
\begin{equation}
\mathbf{B}(f,g)=\left\langle \rho_{1}(f)\mid\rho_{1}(g)\right\rangle -\left\langle \rho_{2}(f)\mid\rho_{2}(g)\right\rangle .\label{eq:BoundaryForm3}
\end{equation}
Hence $\left(\mathbb{C}^{2},\rho_{1},\rho_{2}\right)$ is a boundary
triple for $P.$ The set of selfadjoint restrictions of $P$ is parametrized
by the set of unitary $2\times2$ matrices, see e.g., \cite{dO09}.
Explicitly, any unitary $2\times2$ matrix $B$ determines a selfadjoint
restriction $P_{B}$ of $P$ by setting 
\begin{equation}
\mathscr{D}\left(P_{B}\right):=\left\{ f\in\mathscr{D}\left(P\right)\mid B\rho_{1}(f)=\rho_{2}(f)\right\} .\label{eq:ExtensionDomain1}
\end{equation}
Conversely, every selfadjoint restriction of $P$ is obtained in this
manner. For completeness, we begin by working out the details connecting
the boundary form formulation to the von Neumann deficiency space
approach in our setting. 
\begin{defn}
For $I_{1}$ and $I_{2}$, consider the following two copies of the
two-dimensional Hilbert space $\mathbb{C}^{2}$:
\[
\rho_{2}(\mathscr{D}(P))=\mathscr{B}_{L}=\left\{ \left(\begin{array}{c}
f(0)\\
f(\alpha)
\end{array}\right)\left.\right|f\in\mathscr{D}(P)\right\} ,
\]
and
\[
\rho_{1}(\mathscr{D}(P))=\mathscr{B}_{R}=\left\{ \left(\begin{array}{c}
f(1)\\
f(\beta)
\end{array}\right)\left.\right|f\in\mathscr{D}(P)\right\} ,
\]
where $\rho_{1}$ and $\rho_{2}$ are the respective boundary-restrictions. 
\end{defn}
For $f\in\mathscr{D}(P)$, let $f=\varphi_{0}+f_{+}+f_{-}$ be the
decomposition in (\ref{eq:vN2}), i.e., $\varphi_{0}\in\mathscr{D}_{0}$,
$f_{\pm}\in\mathscr{D}_{\pm}$. Note that, for $f\in\mathscr{D}(P)$,
we have
\begin{equation}
\mathbf{B}(f,f)=\left\Vert f_{R}\right\Vert _{\mathscr{B}_{R}}^{2}-\left\Vert f_{L}\right\Vert _{\mathscr{B}_{L}}^{2}=\left\Vert f_{+}\right\Vert ^{2}-\left\Vert f_{-}\right\Vert ^{2}.\label{eq:bdd-triple-1}
\end{equation}
Set
\[
\mathscr{C}_{b}=\{B:\mathscr{B}_{L}\rightarrow\mathscr{B}_{R},\mbox{ isometric}\},
\]
i.e., a copy of the matrix-group of all unitary $2\times2$ matrices. 
\begin{prop}
\label{prop:vN-triple}In the two-interval example, there is a natural
isomorphism $\mathscr{C}_{b}\cong\mathscr{C}(L)$ where $\mathscr{C}(L)$
is defined as in Lemma \ref{lem:vN def-space}.
\begin{proof}
Note that by Lemma \ref{lem:vN def-space}, each of the selfadjoint
extensions from the minimal domain $\mathscr{D}_{0}$ (zero-boundary
conditions) is determined by some $U\in\mathscr{C}(L)$ via $f_{-}=Uf_{+}$,
$f_{+}\in\mathscr{D}_{+}$. Specializing (\ref{eq:bdd-triple-1})
to the domain $\mathscr{D}_{U}$ in (\ref{eq:vN2}), i.e., $f=\varphi_{0}+f_{+}+Uf_{+}$,
$\varphi_{0}\in\mathscr{D}_{0}$, $f_{+}\in\mathscr{D}_{+}$, we get
\[
\mathbf{B}(f,f)=\left\Vert f_{R}\right\Vert ^{2}-\left\Vert f_{L}\right\Vert ^{2}=\left\Vert f_{+}\right\Vert ^{2}-\left\Vert Uf_{+}\right\Vert ^{2}=0
\]
for all $f\in\mathscr{D}_{U}$. It follows that 
\begin{equation}
f_{L}\mapsto f_{R},
\end{equation}
defined for $f\in\mathscr{D}_{U}$, defines a unique element $B\in\mathscr{C}_{b}$,
and that
\begin{equation}
\mathscr{C}(L)\ni U\mapsto B\in\mathscr{C}_{b}\label{eq:vN-triple}
\end{equation}
is an isomorphism. 
\end{proof}
\end{prop}
Since the two deficiency spaces $\mathscr{D}_{\pm}$ are easy to compute
in each of the two-interval examples, we can write down two isomorphisms
$\mathscr{D}_{+}\overset{b_{+}}{\rightarrow}\mathscr{B}_{L}$, and
$\mathscr{D}_{-}\overset{b_{-}}{\rightarrow}\mathscr{B}_{R}$; note
that all four spaces are two-dimensional. With this, we get the following
representation of the isomorphism in Proposition \ref{prop:vN-triple}:
\begin{cor}[Matrix realization]
\label{cor:matrix-realization}Let $U$, $b_{\pm}$, and $B$ be
the operators in Proposition \ref{prop:vN-triple}, then $B:\mathscr{B}_{L}\rightarrow\mathscr{B}_{R}$
is $B=b_{-}Ub_{+}^{-1}$; see Figure \ref{fig:def-space-isom}.

\begin{center}
\begin{figure}[h]
\begin{centering}
$\xymatrix{\mathscr{D}_{+}\ar[r]^{b_{+}}\ar[d]_{U} & \mathscr{B}_{L}\ar[d]^{B=b_{-}Ub_{+}^{-1}}\\
\mathscr{D}_{-}\ar[r]_{b_{-}} & \mathscr{B}_{R}
}
$
\par\end{centering}

\caption{\label{fig:def-space-isom}}

\end{figure}

\par\end{center}

\end{cor}
With (\ref{eq:vN-triple}) in Proposition \ref{prop:vN-triple}, we
now parametrize the collection of unitary $2\times2$ matrices by
\begin{equation}
B=\left(\begin{array}{cc}
w\: e(\phi) & -\sqrt{1-w^{2}}\: e(\theta-\psi)\\
\sqrt{1-w^{2}}\: e(\psi) & w\: e(\theta-\phi)
\end{array}\right)\label{eq:2-by-2 Unitary}
\end{equation}
where $0\leq w\leq1,$ $\phi,\psi,\theta\in\mathbb{R},$ and 
\begin{equation}
e(x):=e^{i2\pi x}.\label{eq:NormalizedExp}
\end{equation}
Some of our results about the spectrum are summarized in
\begin{thm}
\label{thm:Magic-equation} If $P_{B}=P_{w,\phi,\psi,\theta}$ is
the selfadjoint momentum operator on $[0,1]\cup[\alpha,\beta]$ associated
with the unitary matrix (\ref{eq:2-by-2 Unitary}) via (\ref{eq:ExtensionDomain1}),
then $P_{B}$ has pure point spectrum. In fact, the spectrum of $P_{B}$
is the set of all real solutions $\lambda$ to the equation 
\begin{equation}
\left(e\left(\phi+\lambda\right)-w\right)e\left(\theta-\phi+\left(\beta-\alpha\right)\lambda\right)=w\: e\left(\phi+\lambda\right)-1.\label{eq:Grand-Eigen-Value-Equation}
\end{equation}

\end{thm}
In particular, the spectrum is independent of the parameter $\psi$
and does not depend on the location of the interval $\left[\alpha,\beta\right]$
only on its length. 
\begin{rem}
\textcolor{black}{Our analysis also works when $\alpha=1.$ In this
case one has to be careful with the interpretation of the function
values at the point of intersection of $[0,1]$ and $[\alpha,\beta].$
One way is the interpretations
\[
f(1)=f(1-)=\lim_{t\nearrow1}f(t),\text{ and }f(\alpha)=f(1+)=\lim_{t\searrow1}f(t).
\]
We leave the details to the interested reader. }
\end{rem}

\begin{rem}
For applications, the role of the four parameters in (\ref{eq:2-by-2 Unitary})
is taken up in Section \ref{sec:U(T)} below. Here we merely add a
note about the geometry entailed by the two extreme cases (i) $w=0$,
and (ii) $w=1$. In the first case, the matrix $B$ is off-diagonal,
with zeros in the two diagonal slots, while for (ii) it is diagonal.
As a result, for (ii), we are merely dealing with the orthogonal direct
sum of separate boundary value problems for the two intervals $I_{1}$
and $I_{2}$ treated in isolation; so two separate index $(1,1)$
problems. By contrast, for (i), the two intervals get exchanged each
time a local translation reaches a boundary point. As noted in Section
\ref{sec:U(T)} below, the dynamics in the intermediate case $0<w<1$,
may be computed with the use of a binary random-walk model.
\end{rem}

\section{Spectral Theory\label{sec:Spectral-Theory}}

This section is about a family of selfadjoint operators $P_{B}$ arising
as restrictions to the union of two intervals of the maximal momentum
operator. We explore how these selfadjoint operators $P_{B}$ depend
both on the prescribed boundary matrix $B$ from eq. (\ref{eq:2-by-2 Unitary}),
as well as the two intervals, i.e., on the choice of the two numbers
$\alpha$ and $\beta,$ see eq. (\ref{eq:InnerProduct}). So the operators
and their spectra depend on six parameters in all, and this dependency
is charted in complete detail below.

Note that the difference $\beta\lyxmathsym{\textendash}\alpha$ is
the length of the second interval. Further recall that the unitary
2 by 2 matrix $B$ has four parameters, see (\ref{eq:2-by-2 Unitary});
and that its radial parameter $w$ lies in the closed interval from
0 to 1. We say that the endpoints constitute the extreme cases, and
the open interval $0<w<1,$ the generic case. The role of $w$ is
analyzed further in Section \ref{sec:U(T)} which deals with scattering
theory for the unitary one-parameter groups generated by each of the
operators $P_{B}.$

In the present section, we show that each operator $P_{B}$ has pure
eigenvalue spectrum $\left\{ \lambda_{n}(B)\right\} $ with the eigenvalues
depending on all parameters, i.e., on both $\alpha$ and $\beta,$
and on the matrix $B.$ This dependence is detailed below: In Section
\ref{sub:Boundary-Cases} we deal with two extreme cases for $B;$
Section \ref{sub:Generic-Cases} the generic case; and in Section
\ref{sub:Irrational} we deal with the dichotomy for the value of
the length of the second interval, rational vs irrational. 

We saw in Corollary \ref{cor:matrix-realization} that every selfadjoint
restriction of the maximal momentum operator $P$ has the form $P_{B}$
for some $2\times2$ matrix $B$ as in (\ref{eq:2-by-2 Unitary}).
We begin our investigation of how the spectrum of $P_{B}$ depends
on the parameters $1<\alpha<\beta,$ $\phi,$ $\psi,$ $\theta,$
and $0\leq w\leq1$ by stating two simple lemmas.  

Fix $1<\alpha<\beta.$ Let $I_{1}:=[0,1]$, and $I_{2}:=[\alpha,\beta]$.
Let $e_{x}(y):=e(xy).$ Any complex number $\lambda$ is an eigenvalue
of the maximal momentum operator $P$ in (\ref{eq:MomentumOperator}).
For each complex number $\lambda$ the formal solutions to $Pf=\lambda f$
are 
\begin{equation}
e_{\lambda}^{(a,b)}(x):=\left(a\chi_{I_{1}}(x)+b\chi_{I_{2}}(x)\right)e_{\lambda}(x),\, a,b\in\mathbb{C}.\label{eq:EigenFunction}
\end{equation}

\begin{lem}
\label{lem:The-maximal-operator}The maximal operator $P$ has spectrum
equal to the complex plane and each point in the spectrum is an eigenvalue
with multiplicity two. \end{lem}
\begin{proof}
Since $e_{\lambda}^{(a,b)}$ is in $L^{2}(I_{1}\cup I_{2}),$ it follows
from $Pe_{\lambda}^{(a,b)}=\lambda e_{\lambda}^{(a,b)}$ and (\ref{eq:EigenFunction})
that each $\lambda\in\mathbb{C}$ is an eigenvalue of $P$ with multiplicity
two. \end{proof}
\begin{lem}
The spectrum of any selfadjoint restriction $\widetilde{P}$ of $P$
equals the set of eigenvalues of $\widetilde{P}$ and each eigenvalue
has multiplicity one or two. \end{lem}
\begin{proof}
By Lemma \ref{lem:The-maximal-operator} the spectrum of $\widetilde{P}$
is the $\lambda$ for which $e_{\lambda}^{(a,b)}$ is in the domain
of $\widetilde{P}$ for some choice of constants $a,b.$ The multiplicity
is determined by the dimension of this set of constants. 

Suppose $P_{B}$ is the selfadjoint restriction of $P$ determined
by the matrix (\ref{eq:2-by-2 Unitary}) via (\ref{eq:ExtensionDomain1}).
When the parameters for $B$ are fixed, the expression on the right-hand
side in eq. (\ref{eq:EigenFunction}) is a function of the $x$-variable.
The constants $a$ and $b$ in front of the indicator functions depend
on $B,$ and, in turn, they determine the eigenfunctions, and therefore
the spectrum, of the associated selfadjoint operator $P_{B}.$ 

Suppose $P_{B}$ is the selfadjoint restriction of $P$ determined
by the matrix (\ref{eq:2-by-2 Unitary}) via (\ref{eq:ExtensionDomain1}).
Then the boundary condition $B\rho_{1}(f)=\rho_{2}(f)$ states 
\begin{equation}
\left(\begin{array}{cc}
w\: e(\phi) & -\sqrt{1-w^{2}}\: e(\theta-\psi)\\
\sqrt{1-w^{2}}\: e(\psi) & w\: e(\theta-\phi)
\end{array}\right)\left(\begin{array}{c}
f(1)\\
f(\beta)
\end{array}\right)=\left(\begin{array}{c}
f(0)\\
f(\alpha)
\end{array}\right).\label{eq:ExtensionDomain2}
\end{equation}
Consequently, $\lambda$ is an eigenvalue for $P_{B}$ iff there is
are complex numbers $a,b$ such that $f=e_{\lambda}^{(a,b)}$ satisfies
the boundary condition (\ref{eq:ExtensionDomain2}), i.e., 
\begin{equation}
\left(\begin{array}{cc}
w\: e(\phi) & -\sqrt{1-w^{2}}\: e(\theta-\psi)\\
\sqrt{1-w^{2}}\: e(\psi) & w\: e(\theta-\phi)
\end{array}\right)\left(\begin{array}{c}
a\: e(\lambda)\\
b\: e(\lambda\beta)
\end{array}\right)=\left(\begin{array}{c}
a\\
b\: e(\lambda\alpha)
\end{array}\right).\label{eq:EigenBoundary}
\end{equation}
Since $P_{B}$ is selfadjoint any eigenvalue $\lambda$ must be a
real number. 
\end{proof}
Below we show that the real solutions $\lambda$ to (\ref{eq:EigenBoundary})
is a discrete set having no other accumulation points than plus/minus
infinity.

\subsection{Boundary Cases\label{sub:Boundary-Cases}}

Our analysis depends on the parameter $w.$ We begin by considering
the extreme cases $w=1$ and $w=0.$ 
\begin{thm}[$w=1$]
\label{thm:w=00003D1}Fix $1<\alpha<\beta.$ Let $I_{1}=[0,1]$,
and $I_{2}:=[\alpha,\beta].$ Suppose $P_{B}$ be the selfadjoint
restriction associated with (\ref{eq:2-by-2 Unitary}) via (\ref{eq:ExtensionDomain1})
and $w=1$. Then $P_{B}$ has pure point spectrum $\Lambda_{1}\cup\Lambda_{2}$,
where
\begin{eqnarray}
\Lambda_{1} & = & -\phi+\mathbb{Z}\label{eq:Spectrum (diag) 1}\\
\Lambda_{2} & = & \frac{\phi-\theta}{\beta-\alpha}+\frac{1}{\beta-\alpha}\mathbb{Z}\label{eq:Spectrum (diag) 2}
\end{eqnarray}
and corresponding eigenfunctions are
\begin{equation}
e_{\phi+n}^{\left(1,0\right)}=\chi_{\left[0,1\right]}e_{\phi+n},\;\text{and }e_{\frac{\theta+n}{\beta-\alpha}}^{\left(0,1\right)}=\chi_{\left[\alpha,\beta\right]}e_{\frac{\theta+n}{\beta-\alpha}}.\label{eq:EigenFunction(Diag)}
\end{equation}
Consequently, each eigenvalue has multiplicity one, except the eigenvalues,
if any, in $\Lambda_{1}\cap\Lambda_{2}$ have multiplicity two. \end{thm}
\begin{proof}
In this case the boundary condition (\ref{eq:EigenBoundary}) takes
the form 
\begin{equation}
\left(\begin{array}{cc}
e(\phi) & 0\\
0 & e(\theta-\phi)
\end{array}\right)\left(\begin{array}{c}
a\: e(\lambda)\\
b\: e(\lambda\beta)
\end{array}\right)=\left(\begin{array}{c}
a\\
b\: e(\lambda\alpha)
\end{array}\right).\label{eq:EigenBoundary (diag)}
\end{equation}
That is 
\begin{align*}
a\: e\left(\phi+\lambda\right) & =a\\
b\: e\left(\theta-\phi+\beta\lambda\right) & =b\: e\left(\alpha\lambda\right).
\end{align*}
The stated conclusions are now immediate. \end{proof}
\begin{rem}
It is easy to see that (\ref{eq:Spectrum (diag) 1}) and (\ref{eq:Spectrum (diag) 2})
are the solutions to (\ref{eq:Grand-Eigen-Value-Equation}) when $w=1.$ 
\end{rem}
We present a few examples in order to illustrate some of the possibilities. 
\begin{example}
If $\beta-\alpha=1$ and $\theta=2\phi$ the spectrum equals $-\phi+\mathbb{Z}$
and has uniform multiplicity equal to two. 
\begin{example}
If $\phi=0,\mbox{\ensuremath{\theta}\ is rational}$ and $\beta-\alpha$
is irrational, then the spectrum has uniform multiplicity equal to
one. 
\begin{example}
If $\phi=\theta=0$ and $\beta-\alpha=2,$ then the spectrum is $\frac{1}{2}\mathbb{Z},$
the eigenvalues in $\mathbb{Z}$ have multiplicity two and the non-integer
eigenvalues have multiplicity one. 
\begin{example}
If $\phi=\theta=0$ and $\beta-\alpha$ is irrational, then $0$ has
multiplicity two and the other eigenvalues have multiplicity one. 
\end{example}
\end{example}
\end{example}
In the other extreme case we have \end{example}
\begin{thm}[$w=0$]
\label{thm:w=00003D0}Fix $1<\alpha<\beta.$ Let $I_{1}=[0,1]$,
and $I_{2}:=[\alpha,\beta].$ Suppose $P_{B}$ be the selfadjoint
restriction associated with (\ref{eq:2-by-2 Unitary}) via (\ref{eq:ExtensionDomain1})
and $w=0$. Then $P_{B}$ has pure point spectrum
\begin{equation}
\frac{\frac{1}{2}-\theta}{1+\beta-\alpha}+\frac{1}{1+\beta-\alpha}\mathbb{Z}\label{eq:Spectrum (off-diag)}
\end{equation}
with uniform multiplicity equal to one and an eigenfunction corresponding
to the eigenvalue $\frac{\frac{1}{2}-\theta+n}{1+\beta-\alpha}$,
$n\in\mathbb{Z}$, is 
\begin{equation}
\left(e\left(\frac{1}{2}-\psi+\frac{\left(\frac{1}{2}+n\right)\beta-\left(1-\alpha\right)\theta}{1+\beta-\alpha}\right)\chi_{I_{1}}+\chi_{I_{2}}\right)e_{\frac{\frac{1}{2}-\theta+n}{1+\beta-\alpha}}.\label{eq:EigenFunction (off-diag)}
\end{equation}
\end{thm}
\begin{proof}
In this case the boundary condition (\ref{eq:EigenBoundary}) takes
the form
\begin{equation}
\left(\begin{array}{cc}
0 & -e(\theta-\psi)\\
e(\psi) & 0
\end{array}\right)\left(\begin{array}{c}
a\: e(\lambda)\\
b\: e(\lambda\beta)
\end{array}\right)=\left(\begin{array}{c}
a\\
b\: e(\lambda\alpha)
\end{array}\right).\label{eq:EigenBoundary (off-diag)}
\end{equation}
Writing (\ref{eq:EigenBoundary (off-diag)}) as two equations we get
\begin{align}
-b\: e(\theta-\psi+\lambda\beta) & =a\label{eq:1 (off-diag)}\\
a\: e(\psi+\lambda) & =b\: e(\lambda\alpha).\label{eq:2 (off-diag)}
\end{align}
The first equation shows that $\left|a\right|=\left|b\right|.$ In
particular, we may set $b=1$ and $a=-e(\theta-\psi+\lambda\beta)$.
Then the eigenvalues are determined by the second equation which now
states
\[
-e\left(\theta+\lambda+\beta\lambda\right)=e(\lambda\alpha).
\]
Hence the eigenvalues $\lambda$ are determined by 
\[
\theta+\lambda(1+\beta-\alpha)\in\frac{1}{2}+\mathbb{Z}.
\]
The stated conclusions are now immediate.\end{proof}
\begin{rem}
It is easy to see that (\ref{eq:Spectrum (off-diag)}) are the solutions
to (\ref{eq:Grand-Eigen-Value-Equation}) when $w=0.$ 
\end{rem}

\subsection{Generic Cases\label{sub:Generic-Cases}}

\begin{figure}[H]
\begin{centering}
\begin{tabular}{c}
\includegraphics[scale=0.5]{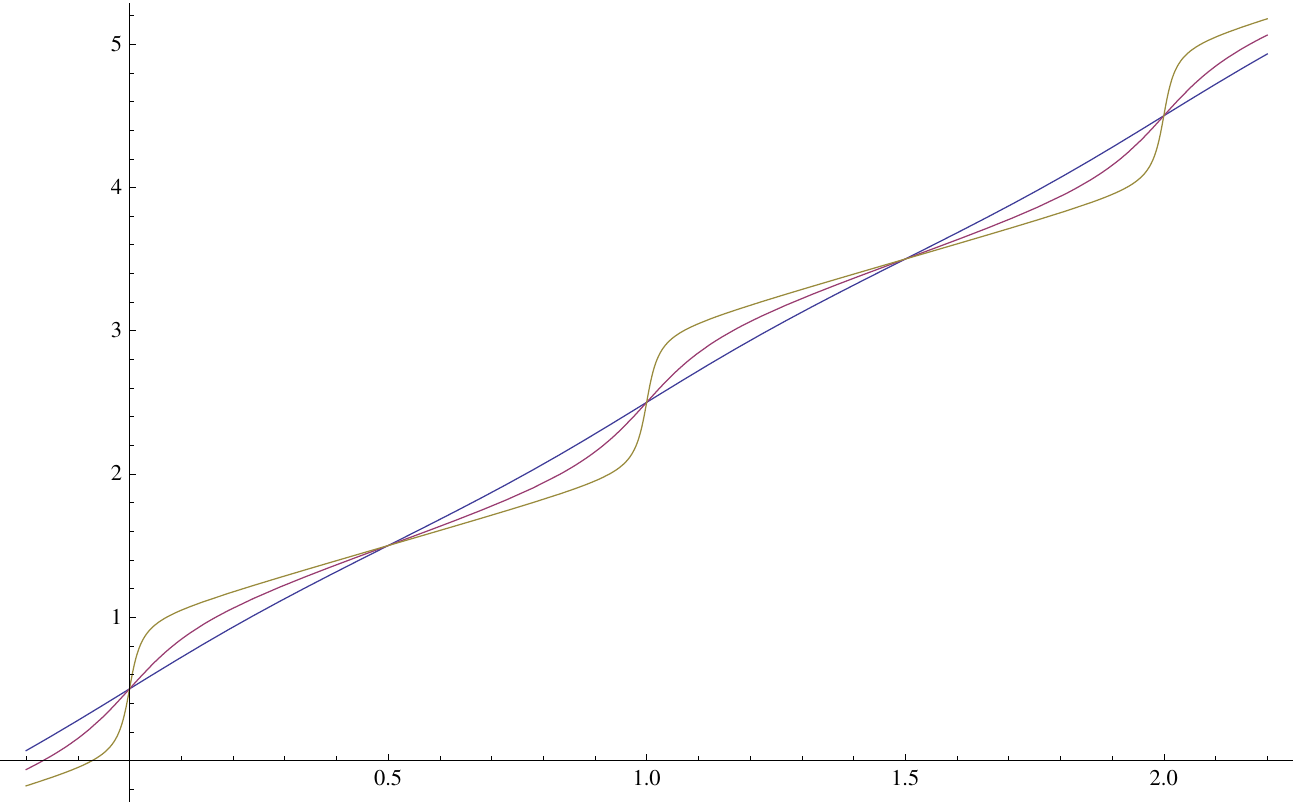}\tabularnewline
\tabularnewline
\end{tabular}
\par\end{centering}

\caption{\label{fig:lambda-vs-w}$y=\theta-\phi+\left(\beta-\alpha\right)t-g(\phi+t)$
for $\theta,\psi,\phi=0,$ $\beta-\alpha=1,$ and $w=0.1,0.5,0.9$. }
\end{figure}

We now consider the non-extreme cases $0<w<1.$ We begin by showing
that the eigenvalues are the solutions to (\ref{eq:Grand-Eigen-Value-Equation}),
thereby completing the proof of Theorem \ref{thm:Magic-equation}. 
\begin{lem}
\label{Lemma:Magic-Equation}Fix $1<\alpha<\beta.$ Let $I_{1}=[0,1]$,
and $I_{2}:=[\alpha,\beta].$ Suppose $P_{B}$ is the selfadjoint
restriction associated with (\ref{eq:2-by-2 Unitary}) via (\ref{eq:ExtensionDomain1})
and $0<w<1$. Then a point $\lambda$ in $\mathbb{R}$ is an eigenvalue
of $P_{B}$ if and only it is a solution to the equation 
\begin{equation}
e\left(\theta-\phi+\left(\beta-\alpha\right)\lambda\right)=\frac{w\: e\left(\phi+\lambda\right)-1}{e\left(\phi+\lambda\right)-w}\label{eq:PhaseIntersection (mix) 1}
\end{equation}
in terms of the argument this is 
\begin{equation}
\theta-\phi+\left(\beta-\alpha\right)\lambda=\frac{1}{2\pi}\arctan\left(\frac{(1-w^{2})\sin\left(2\pi(\phi+\lambda)\right)}{2w-(1+w^{2})\cos\left(2\pi(\phi+\lambda)\right)}\right)+\mathbb{Z}.\label{eq:PhaseIntersection (mix) 2}
\end{equation}
When $\lambda$ is a real solution to (\ref{eq:PhaseIntersection (mix) 1})
then any corresponding eigenfunction is a multiple of 
\begin{equation}
e_{\lambda}^{(a,1)}(x)=\left(a\chi_{I_{1}}(x)+\chi_{I_{2}}(x)\right)e_{\lambda}(x),\label{eq:EigenFunction (mix)}
\end{equation}
where $a$ is determined by (\ref{eq:EigenBoundary3 (mix)}) below.
In particular, the spectrum has uniform multiplicity equal to one. \end{lem}
\begin{proof}
The eigenvalue equation (\ref{eq:EigenBoundary}) is 
\begin{align}
a\: w\: e(\phi+\lambda)-b\sqrt{1-w^{2}}\: e(\theta-\psi+\beta\lambda) & =a\label{eq:EigenBoundary1 (mix)}\\
a\sqrt{1-w^{2}}\: e(\psi+\lambda)+b\: w\: e(\theta-\phi+\beta\lambda) & =b\: e(\alpha\lambda).\label{eq:EigenBoundary2 (mix)}
\end{align}

Equation (\ref{eq:EigenBoundary1 (mix)}) shows $a=0\iff b=0.$ So,
we may set $b=1$, which gives (\ref{eq:EigenFunction (mix)}). 

Solving both (\ref{eq:EigenBoundary1 (mix)}) and (\ref{eq:EigenBoundary2 (mix)})
for $a,$ shows $a$ equals 
\begin{equation}
\frac{\sqrt{1-w^{2}}\: e\left(\theta-\psi+\beta\lambda\right)}{w\: e\left(\phi+\lambda\right)-1}=\frac{1-w\: e\left(\theta-\phi+(\beta-\alpha)\lambda\right)}{\sqrt{1-w^{2}}\: e(\psi+\lambda-\alpha\lambda)}.\label{eq:EigenBoundary3 (mix)}
\end{equation}
Cross-multiplying gives 
\begin{align*}
 & \left(1-w^{2}\right)e\left(\theta+\lambda+\left(\beta-\alpha\right)\lambda\right)\\
 & =\left(1-w\: e\left(\theta-\phi+(\beta-\alpha)\lambda\right)\right)\left(w\: e\left(\phi+\lambda\right)-1\right).
\end{align*}
Expanding the right hand side we see
\begin{align*}
 & \left(1-w^{2}\right)e\left(\theta+\lambda+\left(\beta-\alpha\right)\lambda\right)\\
 & =-1-w^{2}\: e\left(\theta+\lambda+(\beta-\alpha)\lambda\right)+w\:\left(e\left(\theta-\phi+(\beta-\alpha)\lambda\right)+e\left(\phi+\lambda\right)\right).
\end{align*}
Adding $w^{2}\: e\left(\theta+\lambda+\left(\beta-\alpha\right)\lambda\right)$
to both sides gives 
\[
e\left(\theta+\lambda+\left(\beta-\alpha\right)\lambda\right)=-1+w\left(e\left(\theta-\phi+(\beta-\alpha)\lambda\right)+e\left(\phi+\lambda\right)\right).
\]
Subtracting $w\: e\left(\theta-\phi+(\beta-\alpha)\lambda\right)$
yields 
\[
e\left(\theta+\lambda+\left(\beta-\alpha\right)\lambda\right)-w\: e\left(\theta-\phi+(\beta-\alpha)\lambda\right)=-1+w\: e\left(\phi+\lambda\right).
\]
Factoring the left hand side 
\[
\left(e\left(\lambda\right)-w\: e\left(-\phi\right)\right)e\left(\theta+\left(\beta-\alpha\right)\lambda\right)=-1+w\: e\left(\phi+\lambda\right)
\]
or 
\begin{equation}
\left(e\left(\phi+\lambda\right)-w\right)e\left(\theta-\phi+\left(\beta-\alpha\right)\lambda\right)=-1+w\: e\left(\phi+\lambda\right).\label{eq:EigenBoundary4 (mix)}
\end{equation}
Rearranging (\ref{eq:EigenBoundary4 (mix)}) gives (\ref{eq:PhaseIntersection (mix) 1}).
\end{proof}
The following lemma established some simple facts about the right
hand side of (\ref{eq:PhaseIntersection (mix) 1}). 
\begin{lem}
\label{Lemma:Mobius} Fix $0<w<1.$ Consider the M\"{o}bius transformation
\begin{equation}
Mz:=\frac{wz-1}{z-w},\quad z\in\mathbb{C}.\label{eq:Mobius}
\end{equation}
Then $M=M_{w}$ is a bijection of the unit circle $\mathbb{T}$ onto
itself and there is a decreasing continuous bijection $g$ of the
real line $\mathbb{R}$ onto itself such that $g(0)=-1/2$ and $Me(t)=e\left(g\left(t\right)\right)$
for all $t\in\mathbb{R}.$ \end{lem}
\begin{proof}
Recall any M\"{o}bius transformation is a continuous bijection of
the Riemann sphere onto itself. Note $M$ is its own inverse. Taking
the modulus we see 
\[
\left|w\: e\left(t\right)-1\right|=\left|e\left(t\right)-w\right|,\quad t\in\mathbb{R.}
\]
Hence $M$ maps $\mathbb{T}$ into itself and consequently, $M$ is
a bijection of $\mathbb{T}$ onto itself. 

In particular, $\gamma(t):=Me(t),$ $t\in\mathbb{R}$ has minimal
period equal to one and maps any interval of non-zero integer length
onto $\mathbb{T}.$ Since $M$ has positive coefficients $M$ maps
the extended real line onto itself. Since $Mi$ is in the first quadrant
$M$ maps upper half-plane onto itself, in particular, $\gamma\left(\left[0,1/2\right]\right)$
is the part of $\mathbb{T}$ in the upper half-plane. Note $\gamma(0)=M1=-1=e\left(-1/2\right).$ 

By the path lifting theorem there is a unique continuous function
$f:[0,1]\to\mathbb{R}$ such that $f(0)=-1/2$ and $e\left(f(t)\right)=\gamma(t)$
for $t\in[0,1].$ Since $e\left(f(1/2)\right)=\gamma(1/2)=M(-1)=1$
and $e\left(f\left(\left[0,1/2\right]\right)\right)=\gamma\left(\left[0,1/2\right]\right)$
is the part of $\mathbb{T}$ in the upper half-plane, it follows that
$f(1/2)=-1$. Consequently, $f$ is a continuous bijection of $[0,1]$
onto $\left[-\frac{1}{2},-\frac{3}{2}\right].$ 

The desired function is determined by
\[
g(t+n):=f(t)-n
\]
for $0\leq t<1$ and integers $n.$ \end{proof}
\begin{rem}
\textcolor{black}{Since ${\color{black}g(1/4)}$ is in the first quadrant
and ${\color{black}g(1/2)=1}$ then $e\left(g\left[0,1/4\right]\right)$
is the part of $\mathbb{{\color{black}T}}$ in the fourth quadrant
and part of the first quadrant and $e\left(g\left[1/4,1/2\right]\right)$
is the remaining part of $\mathbb{T}$ in the first quadrant. This
``explains'' the general shape of the curve in our pictures below. }
\end{rem}
We are now ready to finish our analysis of the spectrum of $P_{B}.$ 
\begin{thm}[$0<w<1$]
\label{thm:0<w<1}Let the two intervals $I_{1}$ and $I_{2}$ be
as described, with $\alpha$ and $\beta$ fixed and denoting the endpoints
of $I_{2}$, i.e., $1<\alpha<\beta$. Pick parameters $\phi$, $\psi$,
and $w$, $0<w<1,$ describing one of the selfadjoint extensions $P_{B}$.
Then the spectrum of $P_{B}$ is pure point, with the eigenvalue list
equal to the sequence $\{\lambda_{n}\}$ of real numbers solving equation
(\ref{eq:PhaseIntersection (mix) 1}). Every number $\lambda_{n}$
in the spectrum of $P_{B}$ is isolated, and, for $n$ fixed, the
corresponding eigenspace is spanned by the associated eigenfunction
from (\ref{eq:EigenFunction (mix)}).\end{thm}
\begin{proof}
By Lemma \ref{Lemma:Magic-Equation} the eigenvalues of $P_{B}$ are
the real solution to 
\[
e\left(\theta-\phi+\left(\beta-\alpha\right)\lambda\right)=M\: e\left(\phi+\lambda\right)
\]
where $M$ is the M\"{o}bius transformation (\ref{eq:Mobius}). Hence,
if $g$ is the function from the conclusion of Lemma \ref{Lemma:Mobius},
then the eigenvalues are the solution to 
\[
e\left(\theta-\phi+\left(\beta-\alpha\right)\lambda-g(\phi+\lambda)\right)=1.
\]
Hence the eigenvalues are the $\lambda$ for which 
\begin{equation}
\theta-\phi+\left(\beta-\alpha\right)\lambda-g(\phi+\lambda)\in\mathbb{Z}.\label{eq:eigenvalue-line-Mobius}
\end{equation}
Let $h:\mathbb{R}\to\mathbb{R}$ be determined by 
\begin{equation}
h(t):=\theta-\phi+\left(\beta-\alpha\right)t-g(\phi+t).\label{eq:h-def}
\end{equation}
 If $t$ is an integer and $n$ is an integer, then $g(t+n)=g(t)-n,$
consequently, 
\begin{equation}
h\left(\left[t,t+n\right)\right)=\left[h(t),h(t)+\left(1+(\beta-\alpha)\right)n\right).\label{eq:Image of h}
\end{equation}
is an interval of length $n\left(1+\beta-\alpha\right)$. In particular,
it contains at most $1+n\left(1+\beta-\alpha\right)$ integers, hence
$n$ points from the spectrum of $P_{B}.$ 
\end{proof}
Figure \ref{fig:lambda-vs-w} shows three plots of the left hand side
of (\ref{eq:eigenvalue-line-Mobius}) as a function of $t,$ the plots
illustrate that the curve is less linear for larger values of $w.$ 

\textcolor{black}{The spectrum of $P_{B}$ is continuous in $B$ in
the sense that  }
\begin{cor}
\textcolor{red}{\label{cor:Continuity of the spectrum}}\textcolor{black}{For
integers $n,$ let $\lambda_{n}$ denote the solution to $\theta-\phi+\left(\beta-\alpha\right)t-g(\phi+t)=n.$
The the spectrum of $P_{B}$ is the set $\left\{ \lambda_{n}\mid n\in\mathbb{Z}\right\} ,$
$\lambda_{n}<\lambda_{n+1},$ and for fixed $n,$ $(w,\phi,\psi,\theta,\alpha,\beta)\to\lambda_{n}$
is continuous. }The dependence of the eigenvalues on the parameter
$w$ in $B$ is harmonic. Specifically in the formula (\ref{eq:w-possion})
the circle is given by the parameter $t$, and when viewed as a function
on the circle, we therefore get an extension to the inside disk, represented
in polar coordinates $(w,u)$ where $w$ is radius and $u$ is angle;
and with continuity up to $w=1$.
\end{cor}

\begin{figure}[H]
\begin{centering}
\begin{tabular}{c}
\includegraphics[scale=0.5]{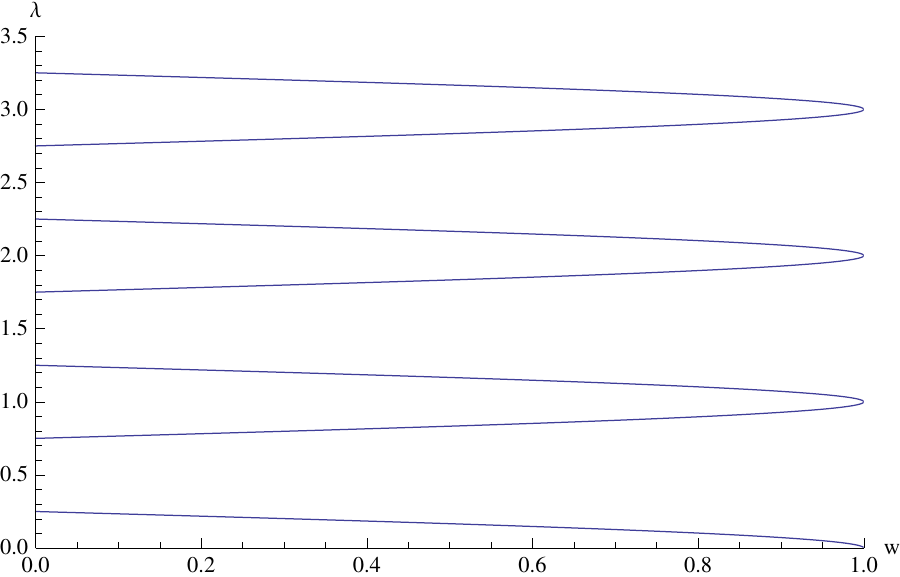}\tabularnewline
\tabularnewline
\end{tabular}
\par\end{centering}

\caption{\label{fig:3.1}$\Omega=(0,1)\cup(2,3)$, and $\theta,\psi,\phi=0$;
$\Lambda_{w=0}=\frac{1}{4}+\frac{1}{2}\mathbb{Z}$, and $\Lambda_{w=1}=\mathbb{Z}$. }
\end{figure}

\begin{figure}[H]
\begin{centering}
\begin{tabular}{c}
\includegraphics[scale=0.5]{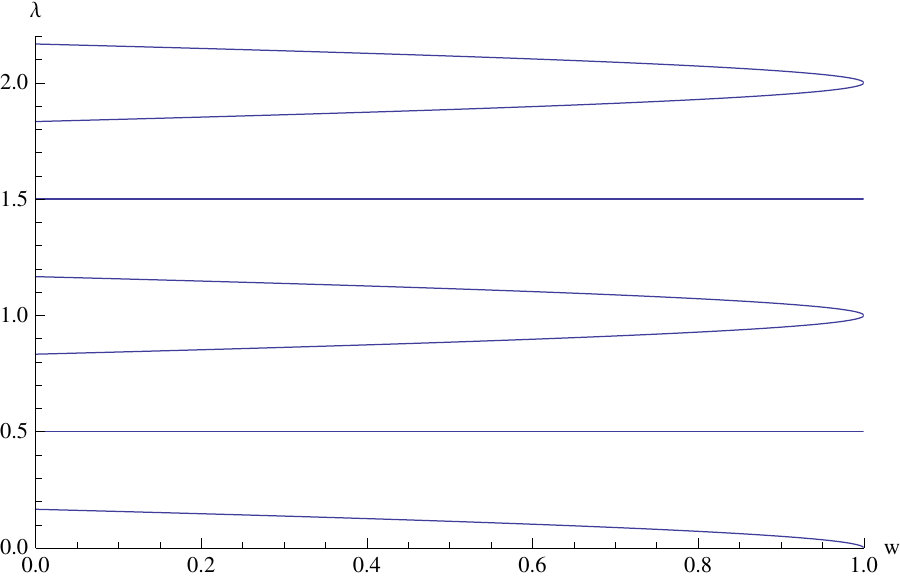}\tabularnewline
\tabularnewline
\end{tabular}
\par\end{centering}

\caption{\label{fig:3.2}$\Omega=(0,1)\cup(2,4)$, and $\theta,\psi,\phi=0$;
$\Lambda_{w=0}=\frac{1}{6}+\frac{1}{3}\mathbb{Z}$, and $\Lambda_{w=1}=\frac{1}{2}\mathbb{Z}$. }
\end{figure}

\begin{figure}[H]
\begin{centering}
\begin{tabular}{c}
\includegraphics[scale=0.5]{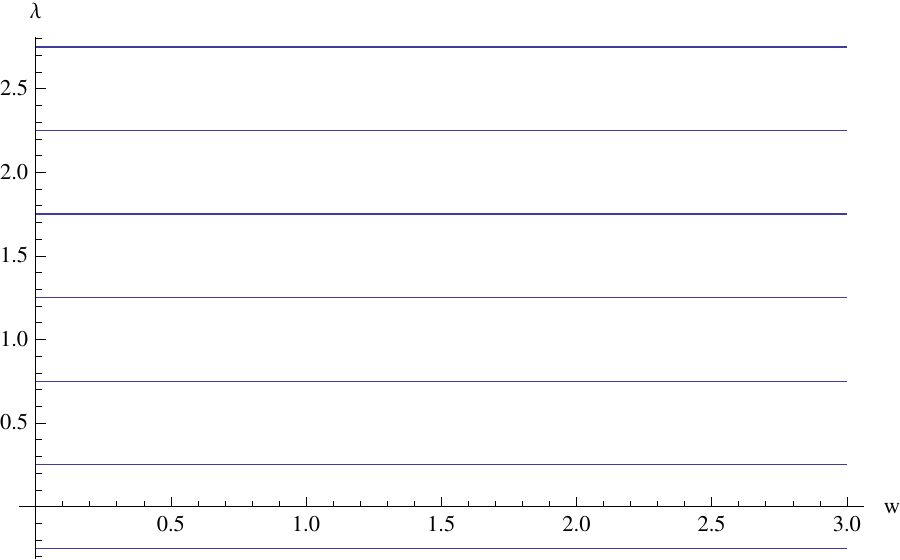}\tabularnewline
\tabularnewline
\end{tabular}
\par\end{centering}

\caption{\label{fig:3.3}$\Omega=(0,1)\cup(2,3)$, $\psi,\phi=0$, and $\theta=1/2$;
$\Lambda_{w}=\frac{1}{2}\mathbb{Z}$, for all $0\leq w\leq1$.}
\end{figure}

\begin{proof}
Let $h$ be as in (\ref{eq:h-def}). Then $\lambda_{n}$ is the first
coordinate of the point of intersection of the curve $y=h(t)$ and
the horizontal line $y=n.$ Clearly the curve is continuous in $\phi,$
$\psi,$ $\theta,$ $\alpha,$ and $\beta.$ By construction, $g$
is the argument of a continuous logarithm. Consequently, 
\begin{equation}
g(t)=-\frac{1}{2}+\frac{1}{2\pi}\mathrm{Im}\int_{0}^{t}\frac{\gamma'}{\gamma}=-\frac{1}{2}-\int_{0}^{t}\frac{1-w^{2}}{1-2w\cos(2\pi u)+w^{2}}du,\label{eq:w-possion}
\end{equation}
where $\gamma(t)=Me(t)$ as in Lemma \ref{Lemma:Mobius}. To get the
harmonic extension, note that the Poisson kernel in (\ref{eq:w-possion})
offers a harmonic extension from the circle to the inside disk. We
use the parameter $t$ for the circle, and the extension is to the
inside disk is in polar coordinates $(w,u)$ with $w$ as radius and
$u$ angle. So the statement about continuity up to $w=1$ is just
continuity up to the boundary circle for the harmonic extension. The
latter is a known property of Poisson kernel.\end{proof}
\begin{rem}
The spectrum of $P_{B}$ is uniform in the sense that for any integer
$n,$ we have 
\[
h\left(\left[\lambda_{n},\lambda_{n}+1\right)\right)=\left[n,n+1+\beta-\alpha\right).
\]
Where $\lambda_{n}$ is as in Corollary \ref{cor:Continuity of the spectrum}.
Consequently, the interval $\left[\lambda_{n},\lambda_{n}+1\right)$
contains exactly $\left\lfloor 1+\beta-\alpha\right\rfloor $ points
from the spectrum. Here, $\left\lfloor 1+\beta-\alpha\right\rfloor $
is an integer such that $\beta-\alpha<\left\lfloor 1+\beta-\alpha\right\rfloor \leq1+\beta-\alpha.$ 
\end{rem}
Points in the spectrum are separated in the sense that 
\begin{cor}
\label{cor:Spectral-separation}There exists $\delta>0,$ so that
if $\lambda\neq\lambda'$ are eigenvalues of $P_{B},$ then $\left|\lambda-\lambda'\right|>\delta.$ \end{cor}
\begin{proof}
The function $t\to\left\langle e_{t}\mid1\right\rangle =\int_{0}^{1}e_{t}+\int_{\alpha}^{\beta}e_{t}$
is continuous and $=1$ at $t=0.$ Hence there is a $\delta>0,$ such
that $|t|\leq\delta$ implies $\left|\left\langle e_{t}\mid1\right\rangle \right|>\frac{1}{2}.$
Since $\lambda-\lambda'$ is a root of $t\to\left\langle e_{t}\mid1\right\rangle $,
we conclude $\left|\lambda-\lambda'\right|>\delta.$
\end{proof}
Let $\lambda_{n}$ be the enumeration of the spectrum of $P$ from
Corollary \ref{cor:Continuity of the spectrum}. The following figures
shows plots of $\lambda_{n}$ as a function of $w$ for several values
of $n.$ Here both the selfadjoint operator $P_{B},$ and the numbers
$\lambda_{n}$ in its eigenvalue list, depend on the choice of the
matrix $B$, but for brevity, in some formulas, we have suppressed
the $B$-dependence. Now $\lambda_{n}=\lambda_{n}(B),$ but in Figure
\ref{fig:3.1}-\ref{fig:3.3} we zoom in on the dependence of the
$w$ parameter from $B.$ Figure \ref{fig:3.1}--\ref{fig:3.3} show
examples where $\delta_{w}\to0$ as $w\to1$ and where $\delta_{w}\geq\widetilde{\delta}>0$
for all $w.$ Here $\delta_{w}=\delta$ is determined in Corollary
\ref{cor:Spectral-separation} above. 

\begin{cor}
\textcolor{black}{Suppose $\beta-\alpha$ is rational and $p,q>0$
are integers such that $\beta-\alpha=p/q.$ Then there is a finite
set $L$ such that the spectrum of $P_{B}$ is $L+q\mathbb{Z}.$ }\end{cor}
\begin{proof}
\textcolor{black}{Suppose $p>0.$ Let $\lambda_{n}$ be as in Corollary
\ref{cor:Continuity of the spectrum}. In particular, $h\left(\lambda_{0}\right)=0$
and 
\[
h\left(\lambda_{0}+q\right)=h\left(\lambda_{0}\right)+\left(1+\beta-\alpha\right)q=q+p.
\]
Hence, setting $t=\lambda_{0}$ and $n=q$ in (\ref{eq:Image of h})
gives 
\[
h\left(\left[\lambda_{0},\lambda_{0}+q\right)\right)=\left[0,q+p\right)
\]
which contains $q+p$ integers. Let $L:=\left\{ \lambda_{0},\lambda_{1},\cdots,\lambda_{q+p-1}\right\} $
be the corresponding eigenvalues of $P_{B}.$ By construction $L\subset\left[\lambda_{0},\lambda_{0}+q\right)$
and the spectrum of $P_{B}$ is $L+q\mathbb{Z}.$ }
\end{proof}

\subsection{Irrational $\beta-\alpha$\label{sub:Irrational}}

In the remainder of this section we discuss the spectrum of $P_{B}$
when $\beta-\alpha$ is irrational. 
\begin{cor}[Asymptotics]
\textcolor{black}{{} }Let $\lambda_{n}$ be the enumeration of the
spectrum of $P$ from Corollary \ref{cor:Continuity of the spectrum}.\textcolor{black}{{}
If $\beta-\alpha$ is irrational, then
\[
\lambda_{0}+\frac{k-1}{1+\beta-\alpha}-1\leq\lambda_{k}<\lambda_{0}+\frac{k-1}{1+\beta-\alpha}+1
\]
for all $k\in\mathbb{Z}.$ }\end{cor}
\begin{proof}
For any real number $x,$ let $\left\lceil x\right\rceil $ be the
real number satisfying $0<\left\lceil x\right\rceil \leq1$ and $x-\left\lceil x\right\rceil \in\mathbb{Z}.$
Well, for any integer $q>0,$ the set $h\left(\left[\lambda_{0},\lambda_{0}+q\right)\right)$
is the interval $[0,\left(1+\beta-\alpha\right)q).$ By irrationality
of $\beta-\alpha,$ the interval $\left[\lambda_{0},\lambda_{0}+q\right)$
contains $\left\lceil \left(1+\beta-\alpha\right)q\right\rceil $
points from the spectrum. Similarly, $\left[\lambda_{0},\lambda_{0}+q+1\right)$
contains $\left\lceil \left(1+\beta-\alpha\right)(q+1)\right\rceil $
points from the spectrum. Consequently,
\[
\lambda_{0}+q\leq\lambda_{\left\lceil \left(1+\beta-\alpha\right)q\right\rceil +1}<\lambda_{\left\lceil \left(1+\beta-\alpha\right)q\right\rceil +2}<\cdots<\lambda_{\left\lceil \left(1+\beta-\alpha\right)(q+1)\right\rceil }<\lambda_{0}+q+1.
\]
Restating we see 
\begin{equation}
\lambda_{0}+q\leq\lambda_{k}<\lambda_{0}+q+1,\label{eq:asymptotics-proof}
\end{equation}
when $\left\lceil \left(1+\beta-\alpha\right)q\right\rceil +1\leq k\leq\left\lceil \left(1+\beta-\alpha\right)(q+1)\right\rceil .$
Now 
\[
\left\lceil \left(1+\beta-\alpha\right)q\right\rceil +1\leq k\leq\left\lceil \left(1+\beta-\alpha\right)(q+1)\right\rceil 
\]
implies 
\[
\left(1+\beta-\alpha\right)q+1\leq k\leq\left(1+\beta-\alpha\right)(q+1)+1.
\]
Solving for $q$ gives
\[
q\leq\frac{k-1}{1+\beta-\alpha}\leq q+1.
\]
The desired conclusions now follow from (\ref{eq:asymptotics-proof}). 

We leave the discussion of the case $q<0$ for the reader. 
\end{proof}
In the general case (\ref{eq:PhaseIntersection (mix) 2}) we study
the intervals in the $\lambda$ axis between neighboring asymptotes:
a sequence of equal-length intervals. The difference $\beta-\alpha$:
rational vs irrational. We now study the points on the axis corresponding
to the intersections in the above examples. By Lemma \ref{Lemma:Magic-Equation}
these points constitute the spectrum $\Lambda$. 
\begin{rem}
\label{Remark: beta-alpha irrational}If \textcolor{black}{$\beta-\alpha$
is irrational, the spectrum is purely aperiodic, see details below
in Example }\ref{Example: dense orbits}.
\end{rem}
The examples below serve to illustrate a number of spectral theoretic
issues accounted for in our theorems.

They are motivated by several considerations:

(i) The examples, both in this section, and in Section \ref{sec:Spectral-Pairs},
show that the possibilities, outlined in the abstract, may indeed
arise.

(ii) A vector bundle over $U(2$): This is illustrated below. Recall,
for each $B$ in $U(2),$ the eigenspaces for $P_{B}$ form a line
bundle over the group $U(2)$ as base-space, with degeneracies at
places when the multiplicity jumps from $1$ to $2;$

(iii) The examples further serve to make concrete the subtle interplay
between geometry and spectrum;

(iv) They illustrate some additional issues dealing with the way a
continuum of parameters has symmetry breaking, and thus produce cases
with properties of special interest; for example illustrating how
the spectral pairs (SP) form a small subfamily; and

\begin{figure}[h]
\begin{centering}
\begin{tabular}{cc}
\includegraphics{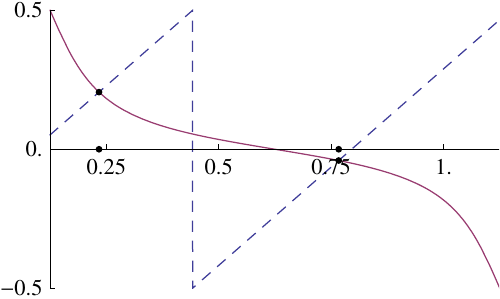} & \includegraphics{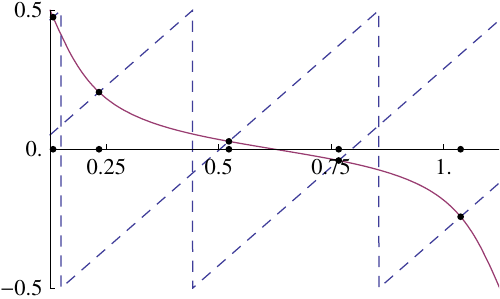}\tabularnewline
$L_{0}$ &  $\cup_{s=0}^{1}(L_{s}-s)$\tabularnewline
\end{tabular}
\par\end{centering}

\caption{\label{fig:dense-orbits} Dense orbits}
\end{figure}

\begin{figure}[h]
\begin{centering}
\begin{tabular}{cc}
\includegraphics{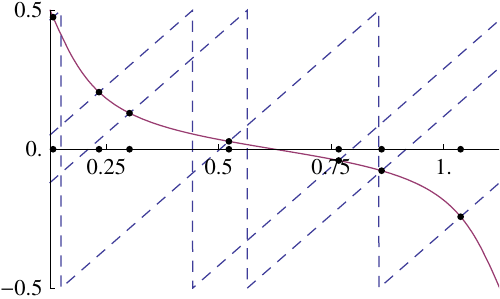} & \includegraphics{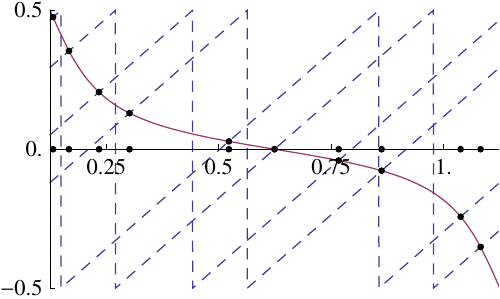}\tabularnewline
 $\cup_{s=0}^{2}(L_{s}-s)$ &  $\cup_{s=0}^{3}(L_{s}-s)$\tabularnewline
\end{tabular}
\par\end{centering}

\caption{\label{fig:dense-orbits-1} Dense orbits}
\end{figure}

\begin{figure}[h]
\begin{centering}
\begin{tabular}{cc}
\includegraphics{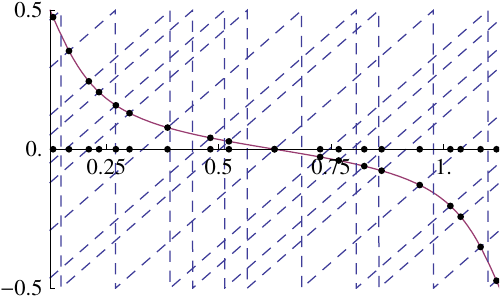} & \includegraphics{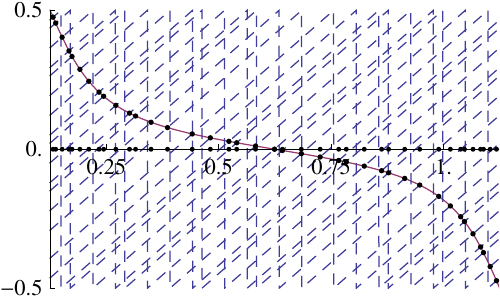}\tabularnewline
 $\cup_{s=0}^{7}(L_{s}-s)$ &  $\cup_{s=0}^{15}(L_{s}-s)$\tabularnewline
\end{tabular}
\par\end{centering}

\caption{\label{fig:dense-orbits-2} Dense orbits: The spectrum as a random-number
generator.}
\end{figure}

(v) The various families of two-interval cases serve as models for
scattering theory in physics (see Section \ref{sub:Applications}2.1),
so they are models for systems with a much richer structure, models
for potential barriers, for quantum jumps, \ldots{} . By model, we
mean, up to unitary equivalence, component-by component.
\begin{example}
\label{Example: dense orbits}Choose $w=\frac{1}{\sqrt{3}}$, $\theta=-1/4$,
$\phi=-1/8$, $\alpha=3$, $\beta=3+\sqrt{2}$. Let $J_{s}=[s+1/8,s+1+1/8)$,
$s\in\mathbb{Z}$, i.e., $J_{s}$ is the $s^{th}$ interval between
two neighboring branch cuts. Let $L_{s}$ be the eigenvalues in $J_{s}$,
then $\Lambda=\bigcup_{s\in\mathbb{Z}}L_{s}$, and the set $\cup_{s\in\mathbb{Z}}\left(L_{s}-s\right)$
is dense in $J_{0}$. See Figures \ref{fig:dense-orbits}-\ref{fig:dense-orbits-2}. 

This conclusion holds whenever $I_{2}$ is chosen with irrational
length $\beta-\alpha$. The reason is that this number then induces
an irrational rotation which is known to have dense orbits; a basic
example in symbolic dynamics. 

Figures \ref{fig:dense-orbits}-\ref{fig:dense-orbits-2} are based
on (\ref{eq:PhaseIntersection (mix) 1}), the argument function from
the proof of Lemma \ref{Lemma:Mobius}, reduction modulo one gives
the branch curves. More precisely, the figures are obtained in a number
of steps. \textbf{Step 1}. Fix $B$ as in (\ref{eq:2-by-2 Unitary}),
and let $P_{B}$ be the corresponding selfadjoint operator. Consider
the curves for argument function from the right hand side of (\ref{Lemma:Magic-Equation})
and the lines with slope $\beta-\alpha$ from the left hand side of
(\ref{Lemma:Magic-Equation}). We are assuming here that this slope
is irrational. \textbf{Step 2}. Identify the asymptotes, and fix an
interval between two branch cuts; for example the interval $J_{0}$
on the axis closest to 0. Note that all the intervals between neighboring
branch cuts have the same fixed unit-length. They extend both to the
left and to the right of $J_{0}$. \textbf{Step 3}. Since the spectrum
of $P_{B}$ is discrete and infinite, it intersects all these intervals
between branch cuts. Now, translate all of these finite intersections
down to $J_{0}$. Conclusion. Since the line-slope is irrational,
the union of all these sets (inside $J_{0}$) is dense in $J_{0}$.
\end{example}

\begin{rem}
When the difference $\beta-\alpha$ is irrational, then the distribution
of points in $\Lambda$ changes as you move around on the axis between
the cut-intervals. See Remark \ref{Remark: beta-alpha irrational},
and Example \ref{Example: dense orbits} above for a detailed discussion.\end{rem}
\begin{cor}
Let $\lambda_{n}$ be the enumeration of the spectrum of $P$ from
Corollary \ref{cor:Continuity of the spectrum}. For a real number
$x,$ let $\left[x\right]$ be the real number satisfying $0\leq\left[x\right]<1.$
If $\beta-\alpha$ is irrational, then $\left[\lambda_{n}\right]\neq\left[\lambda_{m}\right]$
for all $n\neq m$ and the set $\left\{ \left[\lambda_{n}\right]\mid n\in\mathbb{Z}\right\} $
is dense in the interval $[0,1].$ \end{cor}
\begin{proof}
$\lambda_{n}$ is a solution, for $t,$ to 
\begin{equation}
\theta-\phi+(\beta-\alpha)t=n+g(\phi+t)\label{eq:Random1}
\end{equation}
iff $\phi+\lambda_{n}$ is a solution, for $t,$ to 
\[
\theta-(1+\beta-\alpha)\phi+(\beta-\alpha)t=n+g(t).
\]
Consequently, replacing $\phi$ by $0$ corresponds to translating
the spectrum and replacing $\theta$ by $\theta-(1+\beta-\alpha)\phi.$
We will therefore assume $\phi=0$ in (\ref{eq:Random1}). 

Suppose $\left[\lambda_{m}\right]=\left[\lambda_{n}\right],$ Then
$\lambda_{m}=\lambda_{n}+k$ for some integer $k.$ Hence $M\: e\left(\lambda_{m}\right)=M\: e\left(\lambda_{n}\right)$
and therefore Lemma \ref{Lemma:Magic-Equation} implies $e\left(\theta+\left(\beta-\alpha\right)\lambda_{m}\right)=e\left(\theta+\left(\beta-\alpha\right)\lambda_{n}\right).$
Consequently, 
\[
\theta+\left(\beta-\alpha\right)\left(\lambda_{n}+k\right)=\theta+\left(\beta-\alpha\right)\lambda_{n}+j
\]
for some integer $j.$ It follows that $\beta-\alpha=\frac{j}{k}$
. Contradicting that $\beta-\alpha$ is irrational. 

We proceed the proof of the density claim by an alternative description
of the set $\left\{ \left[\lambda_{n}\right]\mid n\in\mathbb{Z}\right\} $.
For real $x$ let $-\frac{1}{2}\leq\left\{ x\right\} <\frac{1}{2}$
be chosen such that $x-\left\{ x\right\} $ is an integer. Let $f$
be the path lifting function from the proof of Lemma \ref{Lemma:Mobius}.
Replacing $t$ in (\ref{eq:Random1}) by $t+m$ where $0\leq t<1,$
we get 
\[
\left\{ \theta+(\beta-\alpha)(t+m)\right\} =1+f(t)
\]
since $\left\{ n+g(t+m)\right\} =\left\{ n+f(t)-m\right\} =1+f(t).$
Let $\ell_{m}(t):=\left\{ \theta+(\beta-\alpha)(t+m)\right\} .$ We
will say $\ell_{m}$ is a line. Then, the set $\left\{ \left[\lambda_{n}\right]\mid n\in\mathbb{Z}\right\} $
is the set of solutions to the equations 
\[
\mbox{\ensuremath{\ell}}_{m}(t)=1+f(t),\quad m\in\mathbb{Z}
\]
in the interval $[0,1).$ 

To show that $\left\{ \left[\lambda_{n}\right]\mid n\in\mathbb{Z}\right\} $
is dense in the interval $[0,1]$ we must show that any subinterval
$[a,b]\subseteq[0,1]$ contain at least one point of the form $\left[\lambda_{n}\right].$
Consider the rectangle $R=[a,b]\times[1+f(b),1+f(a)].$ Note the graph
of $1+f(t),$ $a\leq t\leq b$ is a monotonely decreasing function
passing through the upper left hand corner of $R$ and the lower right
hand corner of $R.$ Since $\left\{ \left[m\left(\beta-\alpha\right)\right]\mid m\in\mathbb{Z}\right\} $
is dense in $[0,1],$ there are lines $\ell_{m}$ passing through
a dense set of point of the bottom $[a,b]\times[1+f(b),1+f(b)]$ of
$R.$ Since the lines $\ell_{m}$ all have fixed positive slope $\beta-\alpha$
we can ensure that they pass through the right hand edge $[b,b]\times[1+f(b),1+f(a)]$
of $R.$ These lines must intersect the graph of $1+f(t),$ the first
coordinates of these intersection points are the points in $\left\{ \left[\lambda_{n}\right]\mid n\in\mathbb{Z}\right\} .$ \end{proof}
\begin{prop}
\label{prop:{finite}+Z}Let $I_{1},I_{2}$ be the two intervals as
before, and $P_{B}$ one of the selfadjoint extensions. Let $\Lambda$
be the spectrum of $P_{B}$. Then $\Lambda=F+\mathbb{Z}$, where $F$
is some finite subset in $\mathbb{R}$, if and only if $\beta-\alpha$
is rational.
\end{prop}

\section{\label{sec:U(T)}Unitary One-Parameter Groups}

Consider the Hilbert space $\mathscr{H}=L^{2}(\Omega)$, $\Omega=I_{1}\cup I_{2}$,
and with the two intervals $I_{1}$ and $I_{2}$ as described. As
before let 
\begin{equation}
B=\left(\begin{array}{cc}
w\: e(\phi) & -\sqrt{1-w^{2}}\: e(\theta-\psi)\\
\sqrt{1-w^{2}}\: e(\psi) & w\: e(\theta-\phi)
\end{array}\right)\label{eq:U(t)-M}
\end{equation}
be a unitary matrix with the parameters $w,\varphi,\theta,\psi$ chosen. 

Consider the Hermitian operator $L=\frac{1}{i2\pi}\frac{d}{dx}$ on
its minimal domain $\mathscr{D}_{0}$ corresponding to function $f\in\mathscr{H}$
such that $f'\in\mathscr{H}$ and $f(0)=f(1)=f(\alpha)=f(\beta)=0$.
Then the selfadjoint extensions $P_{B}$ are prescribed by the choice
in (\ref{eq:U(t)-M}). Specifically, 
\begin{equation}
\mathrm{dom}(P_{B})=\{f\in\mathrm{dom}(P)\left.\right|f_{R}=Bf_{L}\}\label{eq:U(t)-ext-domain}
\end{equation}
where $P$ is the maximal operator and where 
\begin{equation}
f_{R}=\left(\begin{array}{c}
f(1)\\
f(\beta)
\end{array}\right),\;\mbox{and}\; f_{L}=\left(\begin{array}{c}
f(0)\\
f(\alpha)
\end{array}\right).\label{eq:U(t)-bdd-space}
\end{equation}

We get a unitary one-parameter group $U_{B}(t):\mathscr{H}\rightarrow\mathscr{H}$,
$t\in\mathbb{R}$ such that for $f\in\mathrm{dom}(P_{M})$ 
\begin{equation}
\lim_{t\rightarrow0}\frac{1}{t}\left(U_{B}(t)f-f\right)=iP_{B}f.\label{eq:U(t)-iP_M}
\end{equation}

The following is a consequence of the results in Sections \ref{sec:Momentum-Operators}
- \ref{sec:Spectral-Theory}. 
\begin{cor}
\label{cor:U(t)-1}~
\begin{enumerate}
\item Let $f$ be some wave-function localized in $I_{1}$, then if $x$
and $x-t$ are in $I_{1}$, then 
\[
\left(U_{B}(t)f\right)(x)=f(x-t).
\]

\item As the support of $U_{B}(t)f$ hits $x=1$ (the right-hand side boundary
point), then it transfers to the two parts $0$ (left-hand side point
in $I_{1}$), and to $\alpha$ (left-hand side part in $I_{2}$) with
probabilities $w^{2}$ and $1-w^{2}$, respectively. 
\item In the transfer from $x=1$ to $x=0$, the phase is shifted by $e(\varphi)$;
and from $x=1$ to $x=\alpha$, it is phase-shifted by $e(\psi)$. 
\item The boundary conditions (\ref{eq:U(t)-ext-domain}) are preserved
by $U_{B}(t)$ for all $t\in\mathbb{R}$; i.e., we have
\begin{equation}
\left(\begin{array}{c}
\left(U_{B}(t)f\right)(1)\\
\left(U_{B}(t)f\right)(\beta)
\end{array}\right)=B\left(\begin{array}{c}
\left(U_{B}(t)f\right)(0)\\
\left(U_{B}(t)f\right)(\alpha)
\end{array}\right)\label{eq:U(t)-bdd-M}
\end{equation}
for all $f\in\mathrm{dom}(P_{B})$; see also eq. (\ref{eq:U(t)-bdd-space}).\end{enumerate}
\begin{proof}
It follows from (\ref{eq:U(t)-iP_M}) that $\mathrm{dom}(P_{B})$
is invariant under $U_{B}(t)$ for all $t\in\mathbb{R}$, and (\ref{eq:U(t)-bdd-M})
follows from this and Corollary \ref{cor:matrix-realization}.
\end{proof}
\end{cor}
\begin{center}
\begin{figure}[H]
\begin{centering}
\begin{tabular}{c}
\includegraphics{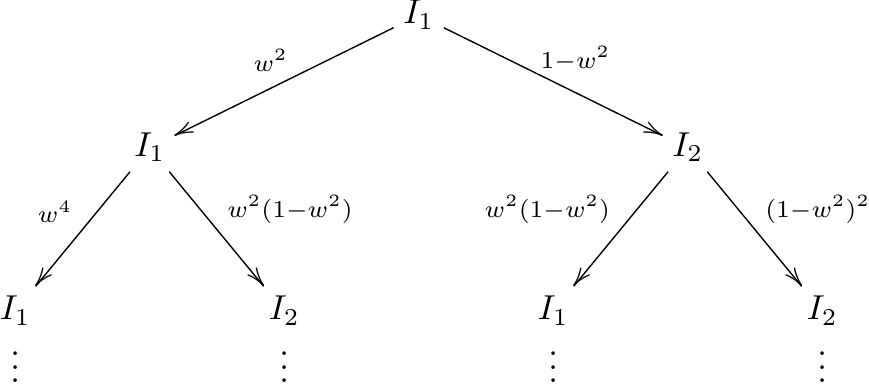}\tabularnewline
\includegraphics{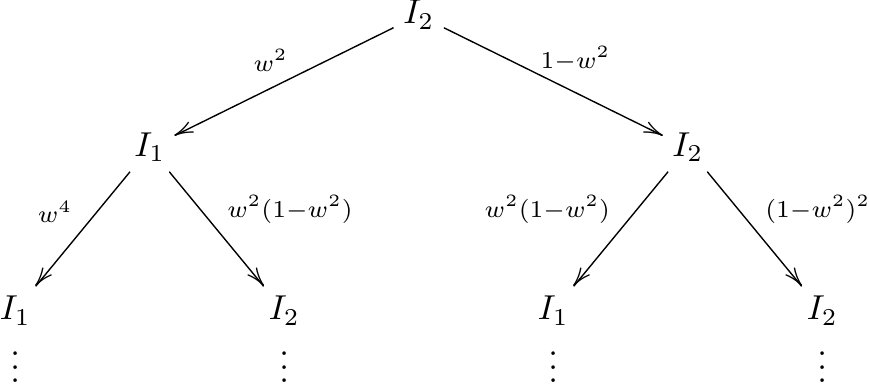}\tabularnewline
\end{tabular}
\par\end{centering}

\caption{Quantum Binary Model}
\end{figure}

\par\end{center}

\section{\label{sec:Spectral-Pairs}Spectral Pairs}

In the study of the spectral pair-problem for subsets in $\mathbb{R}^{d}$
(\cite{Fu74,Jo81,JP98,JP99}) one considers subsets $\Omega$ in $\mathbb{R}^{d}$
of finite positive measure. If the Hilbert space $L^{2}(\Omega)$
has a Fourier basis, we say that $\Omega$ is a \emph{spectral set}.
This means that there is a discrete subset $\Lambda$ in $\mathbb{R}^{d}$
such that $E(\Lambda):=\{e_{\lambda}\left.\right|\lambda\in\Lambda\}$
is an orthogonal basis in $L^{2}(\Omega)$. We then say that $(\Omega,\Lambda)$
is a \emph{spectral pair}, and the set $\Lambda$ is called a \emph{Fourier
spectrum}. The spectrum is never unique. We say $\Omega$ is a \emph{tile,
}or that $\Omega$ tiles by translations, if there is a set $A$ such
that $\sum_{a\in A}\chi_{\Omega}(t-a)=1$ for a.e. $t\in\mathbb{R}^{d}.$
The set $A$ is called a \emph{tiling set} for $\Omega.$ 

An old question from \cite{Fu74} asks if $\Omega$ is spectral if
and only if it tiles $\mathbb{R}^{d}$ by translations. The question
is known to have negative answers for $d\geq3$. But it is open for
$d=2$, and $1$. Another reason for our focus here on $d=1$ and
$\Omega=$ union of two intervals, is that, in this case, we prove
in Corollary \ref{cor:spectral->tile} below that the two properties
are then equivalent. In other words, Fuglede is affirmative in 1D
for these sets.

While the equivalence of the two notions, \textquotedblleft{}tile
by translations\textquotedblright{} and \textquotedblleft{}spectral
set\textquotedblright{} is open for general measurable subsets $\Omega$
in $\mathbb{R}$, of finite positive measure, we show below that,
in case when $\Omega$ is a union of two intervals, then the equivalence
is affirmative. For a discussion of the question for $d\geq3$, see
e.g., \cite{KM06,IKT03,La02,IP98,Ped96,LW97}. The first result proving
non-equivalence was for $d=5$, the paper \cite{Tao04}. Others followed
shortly after, proving non-equivalence for $d=4$, and $d=3$.

Let $1<\alpha<\beta$ and suppose $P_{B}$ is the selfadjoint restriction
of the maximal operator $P$ corresponding to the matrix $B=B(w,\phi,\psi,\theta)$
in (\ref{eq:2-by-2 Unitary}). Let $\Lambda_{B}$ be the spectrum
of $P_{B}.$ By (\ref{eq:EigenFunction}) any eigenfunction of $P_{B}$
is of the form 
\[
\left(a_{\lambda}\chi_{[0,1]}+b_{\lambda}\chi_{[\alpha,\beta]}\right)e_{\lambda}
\]
for some complex numbers $a_{\lambda}$ and $b_{\lambda}.$ We say
$P_{B}$ is \emph{spectral} if $a_{\lambda}=b_{\lambda}$ for all
$\lambda\in\Lambda_{B}.$ 
\begin{lem}
\label{lem:Spectral-Restriction}If $P_{B}$ is spectral, then $\left(\left[0,1\right]\cup\left[\alpha,\beta\right],\Lambda_{B}\right)$
is a spectral pair. Conversely, if $\left(\left[0,1\right]\cup\left[\alpha,\beta\right],\Lambda\right)$
is a spectral pair, then there is a spectral $P_{B}$ such that $\Lambda=\Lambda_{B}.$ \end{lem}
\begin{proof}
Since the eigenfunctions of $P_{B}$ form an orthogonal basis the
first claim is obvious. 

Conversely, suppose $([0,1]\cup[\alpha,\beta],\Lambda)$ is a spectral
pair. Then the closure of 
\[
\sum c_{\lambda}e_{\lambda}\to\sum\lambda c_{\lambda}e_{\lambda}
\]
determines a selfadjoint restriction $\widetilde{P}$ of $P.$ The
spectrum of this selfadjoint restriction is $\Lambda$ and the eigenfunctions
are (the multiples of) the $e_{\lambda}$'s with $\lambda\in\Lambda.$
Since any selfadjoint restriction of $P$ is of the form $P_{B}$
for some $B$ the proof is complete. 
\end{proof}
When $ $$\left(\left[0,1\right]\cup\left[\alpha,\beta\right],\Lambda\right)$
is a spectral pair eq. (\ref{eq:EigenBoundary}) can be solved to
find a corresponding matrix $B.$ 
\begin{thm}
Let $\Lambda_{B}$ be the spectrum of the selfadjoint restriction
$P_{B}$ corresponding to the matrix $B=B(w,\phi,\psi,\theta)$ in
(\ref{eq:2-by-2 Unitary}). If $w=1,$ then $P_{B}$ is not spectral
and consequently $\left([0,1]\cup[\alpha,\beta],\Lambda_{B}\right)$
is not a spectral pair for any choice of the remaining parameters
$\alpha,\beta,\theta,\psi,\phi,w.$ \end{thm}
\begin{proof}
By Theorem \ref{thm:w=00003D1} every eigenfunction has $a_{\lambda}=0$
or $b_{\lambda}=0.$ Hence this is a direct consequence of Lemma \ref{lem:Spectral-Restriction}. 
\end{proof}

\subsection{Two Theorems}

With a choice of two intervals, specified by the endpoints $\alpha$
and $\beta$ for the second interval, we proceed to characterize the
four parameters in the boundary matrix $B$ from (\ref{eq:2-by-2 Unitary})
in Corollary \ref{cor:matrix-realization} which yield a spectral
pair. In other words, we determine when there is an orthogonal Fourier
basis in $L^{2}$ of the two intervals. Our results are Theorems \ref{thm:-spectral-w=00003D0}
and \ref{thm:spectral-0<w<1} below. Together, they exhaust the possibilities
for spectral pairs formed from the union of two intervals.

In particular, we show that all the two-interval spectral pairs fall
in two classes, (i) lattice tilings (Theorems \ref{thm:-spectral-w=00003D0}),
and (ii) non-lattice tilings (Theorems \ref{thm:spectral-0<w<1}).
In both cases, we find the particular instances of spectrum($P_{B}$)
occurring as the second part in a spectral pair. Since we are in one
dimension, (i) represents the cases when the two intervals tile $\mathbb{R}$
by a fixed multiple of the integers $\mathbb{Z}$. In this case, a
Fourier spectrum may be taken to be the dual lattice in $\mathbb{R}$.
We show that, for (ii), the possible Fourier spectra arising as spectrum($P_{B}$)
must be the union of two lattice co-sets.
\begin{thm}
\label{thm:-spectral-w=00003D0}Let $\Lambda_{B}$ be the spectrum
of the selfadjoint restriction $P_{B}$ corresponding to the matrix
$B=B(w,\phi,\psi,\theta)$ in (\ref{eq:2-by-2 Unitary}). If $w=0,$
then $P_{B}$ is spectral iff
\begin{equation}
\frac{\beta}{1+\beta-\alpha}\in\mathbb{N}\label{eq:Spectral-w=00003D0-Domain}
\end{equation}
and the parameters $\phi,\theta$ satisfy 
\begin{equation}
-\psi+\frac{\left(\theta-\frac{1}{2}\right)\left(1-\alpha\right)}{1+\beta-\alpha}\in\mathbb{Z}.\label{eq:Spectral-w=00003D0-B}
\end{equation}
In the affirmative case the spectrum is $\Lambda_{B}=\frac{\frac{1}{2}-\theta}{1+\beta-\alpha}+\frac{1}{1+\beta-\alpha}\mathbb{Z}.$ \end{thm}
\begin{proof}
By Theorem \ref{thm:w=00003D0} and Lemma \ref{lem:Spectral-Restriction}
$P_{B}$ is spectral iff 
\[
-e\left(\theta-\psi+\beta\frac{\frac{1}{2}-\theta+n}{1+\beta-\alpha}\right)=1
\]
 for all $n\in\mathbb{Z}.$ That is iff 
\[
-\psi+\frac{\left(\theta-\frac{1}{2}\right)\left(1-\alpha\right)+\beta\: n}{1+\beta-\alpha}\in\mathbb{Z}
\]
for all $n\in\mathbb{Z}.$ Considering $n=0$ and $n=1$ it follows
that $P_{B}$ is spectral iff the conditions (\ref{eq:Spectral-w=00003D0-B})
are satisfied. 

The formula for the spectrum is from Theorem \ref{thm:w=00003D0}. \end{proof}
\begin{cor}
If (\ref{eq:Spectral-w=00003D0-Domain}) and (\ref{eq:Spectral-w=00003D0-B})
then $\left([0,1]\cup[\alpha,\beta],\frac{1}{1+\beta-\alpha}\mathbb{Z}\right)$
is a spectral pair and $[0,1]\cup[\alpha,\beta]$ is a tile with tiling
set $(1+\beta-\alpha)\mathbb{Z}.$ \end{cor}
\begin{proof}
If $k=\frac{\beta}{1+\beta-\alpha}$, then 
\[
\left[\alpha-(1+\beta-\alpha)(k-1),\beta-(1+\beta-\alpha)(k-1)\right]=[1,1+\beta-\alpha],
\]
consequently, $[0,1]\cup[\alpha,\beta]$ is a tile with tiling set
$(1+\beta-\alpha)\mathbb{Z},$ when $k$ is an integer. \end{proof}
\begin{rem}
In Theorem \ref{thm:-spectral-w=00003D0} the condition (\ref{eq:Spectral-w=00003D0-Domain})
restricts the geometry of the intervals. Specific examples of sets
satisfying this condition includes $[0,1]\cup[\frac{5}{2},3]$ and
$[0,1]\cup\left[7\pi+1,14\pi\right].$ The other condition in (\ref{eq:Spectral-w=00003D0-B})
is a restriction on the set of selfadjoint restrictions, that is on
the values of $\theta$ and $\psi,$ that determine $B$ and the spectrum. 
\end{rem}
For every $\beta-\alpha\in\mathbb{Z}_{+}$, we may construct tiling
spectral pairs. In fact, fix $p=\beta-\alpha$, and set $\beta=2(1+p)$,
$\alpha=2+p$. Then 
\[
\Omega=[0,1]\cup[2+p,2p+2]
\]
tiles with $(1+p)\mathbb{Z}$. 
\begin{example}
Suppose $\beta-\alpha=2$. The condition $\frac{\beta}{1+\beta-\alpha}=\frac{\beta}{3}\in\mathbb{Z}$,
see (\ref{eq:Spectral-w=00003D0-B}), implies $\beta\in3\mathbb{Z}$.
For $\beta=6$, the set $[0,1]\cup[4,6]$ tiles $\mathbb{R}$ with
$3\mathbb{Z}$; see Figure \ref{fig:3Z tile}. By contrast, $[0,1]\cup[2,4]$
does not tile, and note that (\ref{eq:Spectral-w=00003D0-B}) does
not hold.

\begin{figure}[H]
\begin{centering}
\begin{tabular}{c}
\includegraphics{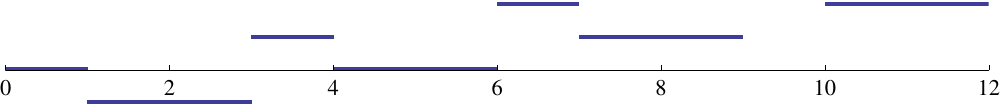}\tabularnewline
\tabularnewline
\end{tabular}
\par\end{centering}

\caption{\label{fig:3Z tile}$\Omega\cup(\Omega+3)\cup(\Omega-3)\cup(\Omega-6)$,
where $\Omega=[0,1]\cup[4,6]$. }
\end{figure}

\end{example}

\begin{example}
Suppose $\beta-\alpha=3$. Then $\frac{\beta}{1+\beta-\alpha}=\frac{\beta}{4}\in\mathbb{Z}$
implies $\beta\in4\mathbb{Z}$. For $\beta=8$, the set $[0,1]\cup[5,8]$
tiles with $4\mathbb{Z}$, by Theorem \ref{thm:-spectral-w=00003D0}.

\begin{figure}[H]
\begin{centering}
\begin{tabular}{c}
\includegraphics{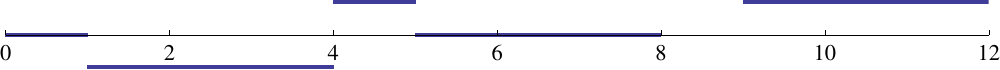}\tabularnewline
\tabularnewline
\end{tabular}
\par\end{centering}

\caption{\label{fig:4Z tile}$\Omega\cup(\Omega+4)\cup(\Omega-4)$, where $\Omega=[0,1]\cup[5,8]$. }
\end{figure}

\end{example}
\begin{thm}
\label{thm:spectral-0<w<1}Let $\Lambda_{B}$ be the spectrum of the
selfadjoint restriction $P_{B}$ corresponding to the matrix $B=B(w,\phi,\psi,\theta)$.
If $0<w<1,$ then $P_{B}$ is spectral iff
\begin{equation}
1<\alpha\text{ is an integer, }\beta=\alpha+1,\label{eq:Spectral-0<w<1-domain}
\end{equation}
 and the parameters $\theta,\phi,\psi,$ and $w$ satisfies
\begin{equation}
\theta-2\phi\in\mathbb{Z},\psi+\left(\alpha-1\right)\phi\in\frac{1}{2}\mathbb{Z},\text{ and }w\in\left\{ \left.\cos\left(2\pi\frac{1+2k}{4\alpha}\right)\right|k\in\mathbb{Z}\right\} .\label{eq:Spectral-0<w<1-B}
\end{equation}
In the affirmative case $\Lambda_{B}=\left\{ -\phi\pm\frac{\cos^{-1}(w)}{2\pi}\right\} +\mathbb{Z},$
hence $\Lambda_{B}=\left\{ -\phi\pm\frac{1+2k}{4\alpha}\right\} +\mathbb{Z},$
for some integer $k.$ \end{thm}
\begin{proof}
We saw earlier that $\lambda$ is a eigenvalue for $P_{B}$ iff 
\begin{equation}
e\left(\theta-\phi+\left(\beta-\alpha\right)\lambda\right)=\frac{w\: e\left(\phi+\lambda\right)-1}{e\left(\phi+\lambda\right)-w}\label{eq:EV}
\end{equation}
and if $\lambda$ is an eigenvalue, then $a_{\lambda}=b_{\lambda}=1$
iff 
\begin{equation}
\sqrt{1-w^{2}}\: e\left(\psi+\lambda-\alpha\lambda\right)=1-w\: e\left(\theta-\phi+\left(\beta-\alpha\right)\lambda\right).\label{eq:a=00003Db=00003D1}
\end{equation}
Taking the square of the modulus to both sides of (\ref{eq:a=00003Db=00003D1})
we get 
\[
0=2w\left(w-\mathrm{c}\left(\theta-\phi+\left(\beta-\alpha\right)\lambda\right)\right).
\]
Where $\mathrm{c}(t):=\cos(2\pi t).$ Since $w\neq0,$ we conclude
\begin{equation}
w=\mathrm{c}\left(\theta-\phi+\left(\beta-\alpha\right)\lambda\right).\label{eq:c=00003Dw}
\end{equation}
Replacing $w$ in (\ref{eq:EV}) by $\mathrm{c}\left(\theta-\phi+\left(\beta-\alpha\right)\lambda\right)$
and simplifying using $\mathrm{s}(t):=\sin(2\pi t)$ and $e(t)=\mathrm{c}(t)+i\mathrm{s}(t)$
we arrive at 
\begin{align*}
 & i\left(\mathrm{c}\left(\phi+\lambda\right)-\mathrm{c}\left(\theta-\phi+\left(\beta-\alpha\right)\lambda\right)\right)\mathrm{s}\left(\theta-\phi+\left(\beta-\alpha\right)\lambda\right)\\
 & =-1+\mathrm{c}^{2}\left(\theta-\phi+\left(\beta-\alpha\right)\lambda\right)+\mathrm{s}\left(\phi+\lambda\right)\mathrm{s}\left(\theta-\phi+\left(\beta-\alpha\right)\lambda\right).
\end{align*}
Consequently, 
\begin{align}
\left(\mathrm{c}\left(\phi+\lambda\right)-\mathrm{c}\left(\theta-\phi+\left(\beta-\alpha\right)\lambda\right)\right)\mathrm{s}\left(\theta-\phi+\left(\beta-\alpha\right)\lambda\right) & =0\label{eq:real}\\
\mathrm{c}^{2}\left(\theta-\phi+\left(\beta-\alpha\right)\lambda\right)+\mathrm{s}\left(\phi+\lambda\right)\mathrm{s}\left(\theta-\phi+\left(\beta-\alpha\right)\lambda\right) & =1.\label{eq:imaginary}
\end{align}
Since $0<w<1,$ if follows from (\ref{eq:c=00003Dw}) that $s\left(\theta-\phi+\left(\beta-\alpha\right)\lambda\right)\neq0,$
hence (\ref{eq:real}) implies 
\[
\mathrm{c}\left(\phi+\lambda\right)=\mathrm{c}\left(\theta-\phi+\left(\beta-\alpha\right)\lambda\right)=w.
\]
It follows from (\ref{eq:c=00003Dw}) and (\ref{eq:imaginary}) that
\[
\mathrm{s}\left(\phi+\lambda\right)=\mathrm{s}\left(\theta-\phi+\left(\beta-\alpha\right)\lambda\right)=\pm\sqrt{1-w^{2}}.
\]
 Consequently, 
\begin{equation}
e\left(\theta-\phi+\left(\beta-\alpha\right)\lambda\right)=e\left(\phi+\lambda\right)=w\pm i\sqrt{1-w^{2}}.\label{eq:spectral-magic}
\end{equation}
is one of the two fixed points $w\pm i\sqrt{1-w^{2}}$ of the M\"{o}bius
transformation $M=M_{w}$ from Lemma \ref{Lemma:Mobius}. Note (\ref{eq:spectral-magic})
implies (\ref{eq:EV}). Hence, if $\lambda$ satisfies (\ref{eq:spectral-magic}),
then $\lambda$ is in $\Lambda_{B}.$ Plugging (\ref{eq:spectral-magic})
into (\ref{eq:a=00003Db=00003D1}) we get 
\begin{equation}
e\left(\psi+\lambda-\alpha\lambda\right)=\sqrt{1-w^{2}}\pm iw=\pm i\left(w\mp i\sqrt{1-w^{2}}\right).\label{eq:a=00003Db=00003D1-ii}
\end{equation}
We have seen that $\lambda$ satisfies (\ref{eq:EV}) and (\ref{eq:a=00003Db=00003D1})
iff $\lambda$ satisfies (\ref{eq:spectral-magic}) and (\ref{eq:a=00003Db=00003D1-ii}). 

Let $\widetilde{\lambda}:=\cos^{-1}(w)/2\pi.$ Then 
\[
e\left(\pm\widetilde{\lambda}\right)=w\pm i\sqrt{1-w^{2}}
\]
 by (\ref{eq:spectral-magic}). And $\mathrm{c}\left(\phi+\lambda\right)=w,$
implies $\lambda\in\left\{ -\phi\pm\widetilde{\lambda}\right\} +\mathbb{Z},$
hence 
\begin{equation}
\Lambda_{B}\subseteq\left\{ -\phi\pm\widetilde{\lambda}\right\} +\mathbb{Z}.\label{eq:spectral-inclusion}
\end{equation}
By \cite{Landau} $\Lambda_{B}$ has density $1+\beta-\alpha>1.$
Hence one of 
\[
\Lambda_{B}\cap\left(\left\{ -\phi+\widetilde{\lambda}\right\} +\mathbb{Z}\right)\text{ and }\Lambda_{B}\cap\left(\left\{ -\phi-\widetilde{\lambda}\right\} +\mathbb{Z}\right)
\]
 has density $>1/2.$ Assume it is $\Lambda_{B}\cap\left(\left\{ -\phi+\widetilde{\lambda}\right\} +\mathbb{Z}\right).$ 

Let $Z_{+}$ be the integers $k$ such that $\lambda_{k}^{+}:=-\phi+\widetilde{\lambda}+k\in\Lambda_{B}$
and similarly, let $Z_{-}$ be the integers $k$ such that $\lambda_{k}^{-}:=-\phi-\widetilde{\lambda}+k\in\Lambda_{B}.$
Plugging $\lambda=\lambda_{k}^{+}=-\phi+\widetilde{\lambda}+k$ into
(\ref{eq:a=00003Db=00003D1-ii}) gives 
\[
e\left(\psi-\left(\phi+\widetilde{\lambda}\right)\left(1-\alpha\right)+\left(1-\alpha\right)k\right)=ie\left(-\widetilde{\lambda}\right)=e\left(\frac{1}{4}-\widetilde{\lambda}\right).
\]
for $k\in Z_{+}.$ Rewriting we get 
\[
e\left(\left(1-\alpha\right)k\right)=e\left(\frac{1}{4}+\alpha\widetilde{\lambda}\right)e\left(\phi\left(1-\alpha\right)-\psi\right)
\]
for $k\in Z_{+}.$ Since the right hand side is independent of $k\in Z_{+}$
and $Z_{+}$ has density $>1/2,$ $1-\alpha$ is an integer. Similarly,
it follows from (\ref{eq:spectral-magic}) that $\beta-\alpha$ is
an integer. The set $\left\{ -\phi\pm\widetilde{\lambda}\right\} +\mathbb{Z}$
has density equal to two, hence the subset $\Lambda_{B}$ has density
at most two, so by \cite{Landau} $1+\beta-\alpha\leq2.$ It follows
that $\beta=1+\alpha$ and $\Lambda_{B}=\left\{ -\phi\pm\widetilde{\lambda}\right\} +\mathbb{Z}.$
In particular, $Z_{+}=Z_{-}=\mathbb{Z},$ $\alpha$ is an integer
$>1$ and $\beta=1+\alpha.$ 

Using $\beta-\alpha=1$ and $\lambda=-\phi\pm\widetilde{\lambda}+k$
for $k\in\mathbb{Z}$ we can write (\ref{eq:spectral-magic}) as 
\[
e\left(\theta-2\phi\pm\widetilde{\lambda}\right)=e\left(\pm\widetilde{\lambda}\right),
\]
thus $\theta-2\phi$ is an integer. Similarly, (\ref{eq:a=00003Db=00003D1-ii})
becomes
\[
e\left(\psi+\left(\alpha-1\right)\phi\right)=e\left(\pm\frac{1}{4}\pm\alpha\widetilde{\lambda}\right).
\]
Consequently, 
\[
\pm\frac{1}{4}\pm\alpha\widetilde{\lambda}\in\psi+\left(\alpha-1\right)\phi+\mathbb{Z}.
\]
Adding these equations shows $2\left(\psi+\left(\alpha-1\right)\phi\right)\in\mathbb{Z}.$
Subtracting them gives $2\alpha\widetilde{\lambda}\in-\frac{1}{2}+\mathbb{Z}.$
But this is equivalent to
\[
\alpha\widetilde{\lambda}\in=\left\{ \frac{-1+2k}{4}\mid k\in\mathbb{Z}\right\} .
\]
This completes the proof. \end{proof}
\begin{cor}
If (\ref{eq:Spectral-0<w<1-domain}) and (\ref{eq:Spectral-0<w<1-B}),
then $\left([0,1]\cup[\alpha,\beta],\left\{ \pm\frac{\cos^{-1}(w)}{2\pi}\right\} +\mathbb{Z}\right)$
is a spectral pair and $[0,1]\cup[\alpha,\beta]$ is a tile with tiling
set $\left\{ 0,1,2,\ldots,\alpha-1\right\} +2\alpha\mathbb{Z}.$ \end{cor}
\begin{rem}
In Theorem \ref{thm:spectral-0<w<1} the condition (\ref{eq:Spectral-0<w<1-domain})
restricts the geometry of the intervals. The set $[0,1]\cup[2,3]$
satisfies (\ref{eq:Spectral-0<w<1-domain}) but not (\ref{eq:Spectral-w=00003D0-Domain}).
The set $[0,1]\cup\left[7\pi+1,14\pi\right]$ satisfies (\ref{eq:Spectral-w=00003D0-Domain})
but not (\ref{eq:Spectral-0<w<1-domain}). Finally, $[0,1]\cup[3,4]$
satisfies both (\ref{eq:Spectral-w=00003D0-Domain}) and (\ref{eq:Spectral-0<w<1-domain}).
The other condition in Theorem \ref{thm:spectral-0<w<1} (\ref{eq:Spectral-0<w<1-B})
is a restriction on the set of selfadjoint restrictions, that is on
the values of $w,\theta,\phi$ and $\psi,$ that determine $B$ and
the spectrum. 
\begin{rem}
(i) Let $\Omega:=[0,1]\cup[\alpha,\beta].$ Suppose (\ref{eq:Spectral-w=00003D0-Domain})
or (\ref{eq:Spectral-0<w<1-domain}) is satisfied, then $\Omega$
is a spectral set. Any choice of a unitary matrix gives us a spectrum
$\Lambda_{B}.$ Most of these spectra will not satisfy (\ref{eq:Spectral-w=00003D0-B})
or (\ref{eq:Spectral-0<w<1-B}). Hence for most $B$ the pair $(\Omega,\Lambda_{B})$
is not a spectral pair, even when $\Omega$ is a spectral set. 

(ii) Suppose $\left(\Omega,\Lambda\right)$ is a spectral pair. By
Lemma \ref{lem:Spectral-Restriction} there is a unitary matrix $B,$
such that $P_{B}$ is spectral and $\Lambda_{B}=\Lambda.$ Construct
a new unitary matrix $B_{0}$ by replacing $\psi$ by a $\psi_{0},$
(keeping $w,\phi,$ and $\theta$ fixed) such that the appropriate
one of (\ref{eq:Spectral-w=00003D0-B}) or (\ref{eq:Spectral-0<w<1-B})
fail. Then $P_{B_{0}}$ is not spectral, yet, since the spectrum is
independent of $\psi$ we have $\Lambda_{B_{0}}=\Lambda_{B}=\Lambda.$ 

(iii) If $\left(\Omega,\Lambda\right)$ is a spectral pair we can
choose selfadjoint restrictions $P_{1}$ and $P_{2}$ of the maximal
operator $P$ such that $\Lambda=\Lambda_{P_{1}}=\Lambda_{P_{2}}$
and the functions $e_{\lambda},\lambda\in\Lambda$ are eigenfunctions
for $P_{1}$ but not for $P_{2}.$  
\end{rem}
\end{rem}
\begin{cor}
\label{cor:Spectral-Set-Criterion} The set $[0,1]\cup[\alpha,\beta]$
is a spectral set iff either $\frac{\beta}{1+\beta-\alpha}$ is an
integer, or $\alpha>1$ is an integer and $\beta=\alpha+1.$ 
\begin{cor}
\label{cor:spectral->tile}If the set $[0,1]\cup[\alpha,\beta]$ is
a spectral set, then it tiles the real line by translations.
\begin{cor}
If $[0,1]\cup[\alpha,\beta]$ is a spectral set and $\beta-\alpha>1,$
then $[0,1]\cup[\alpha,\beta]$ is a lattice tile. 
\end{cor}
\end{cor}
\end{cor}

\subsection{Applications}

While our results for spectral pairs built with two intervals offer
the list of possibilities in terms of the two parameters $\alpha$,
$\beta$, and the boundary matrix (\ref{eq:2-by-2 Unitary}), it is
helpful to visualize some of the possibilities, and discuss geometric
implications. For example, even though the sets spectrum($P_{B}$)
is determined in all cases, its properties are reflected very concretely
by curve intersections, see Figures \ref{fig:Class-B} through \ref{fig:(0,1)(3,6)}
below. And our brief discussions of the examples below serve to highlight
main points. For example, one arrives at a visual for the fact that
the difference $\beta-\alpha$ represents a slope. Also when it is
irrational the restriction to spectrum($P_{B}$) of natural homomorphism
from $\mathbb{R}$ to the quotient $\mathbb{R}/\mathbb{Z}$ (the circle
group) has dense range. See Figures \ref{fig:rat v.s. irrat 2} and
\ref{fig:dense-orbits}.

\begin{figure}[h]
\begin{centering}
\begin{tabular}{cc}
\includegraphics{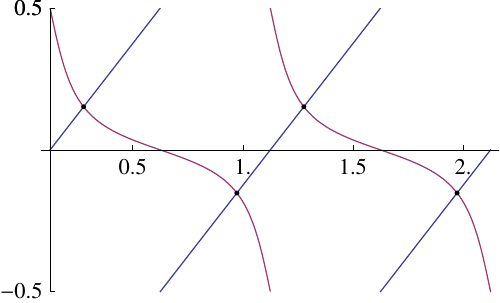} & \includegraphics{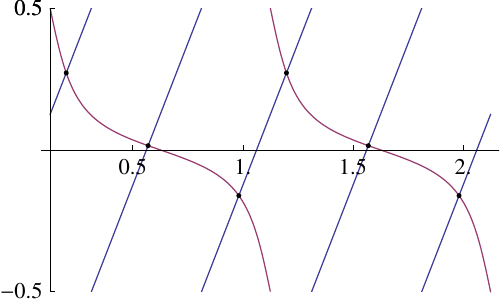}\tabularnewline
 $\beta=4$ &  $\beta=5$\tabularnewline
 & \tabularnewline
\includegraphics{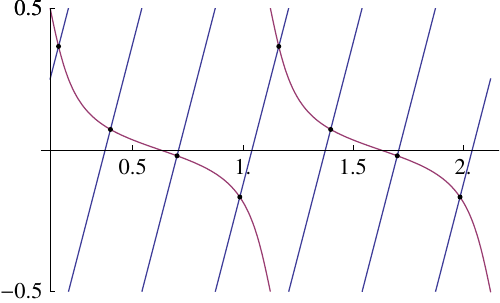} & \includegraphics{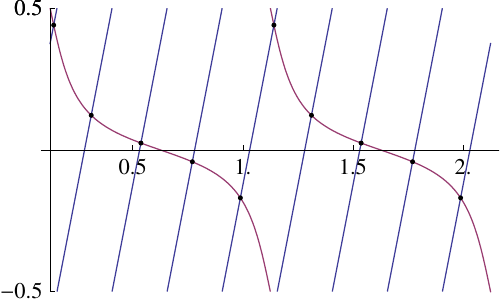}\tabularnewline
 $\beta=6$ & $\beta=7$\tabularnewline
\end{tabular}
\par\end{centering}

\caption{\label{fig:rat v.s. irrat 1}$\beta-\alpha$ is rational, and the
spectrum is periodic.}
\end{figure}

\begin{figure}[h]
\centering{}%
\begin{tabular}{cc}
\includegraphics{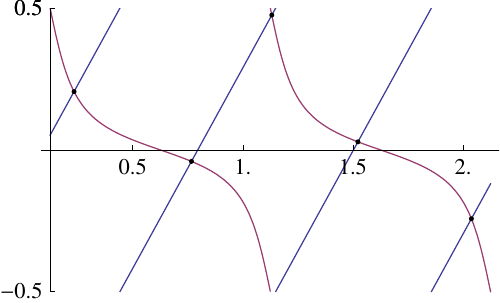} & \includegraphics{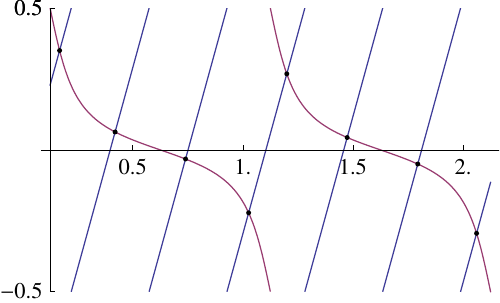}\tabularnewline
$\beta=3+\sqrt{2}$ & $\beta=3+2\sqrt{2}$\tabularnewline
 & \tabularnewline
\includegraphics{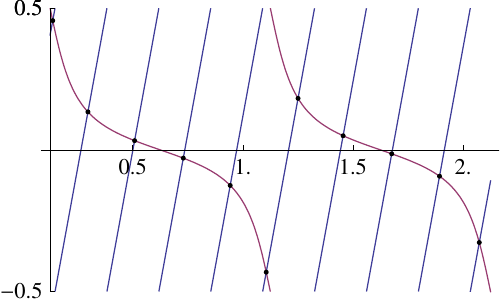} & \includegraphics{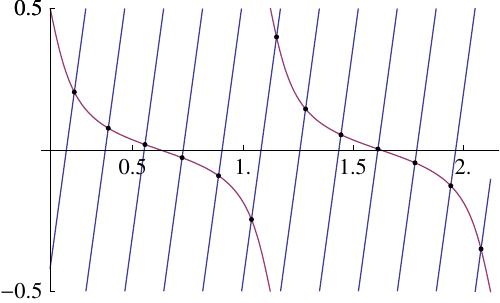}\tabularnewline
$\beta=3+3\sqrt{2}$ & $\beta=3+4\sqrt{2}$\tabularnewline
\end{tabular}\caption{\label{fig:rat v.s. irrat 2}$\beta-\alpha$ is rational, and the
spectrum is purely aperiodic.}
\end{figure}

In the discussion below, we further illustrate how the possibility
for spectral type for our 6-parameter family of selfadjoint operators
depend the roots of modulus $1$ in certain characteristic polynomials;
see Example \ref{example:(0,1)(2,4)}, Figure \ref{fig:(0,1)(2,4)};
Example \ref{example:(0,1)(2,5)}, Figure \ref{fig:(0,1)(2,5)}; and
Example \ref{example:(0,1)(3,5)}, Figure \ref{fig:(0,1)(3,5)}. In
detail; depending on the relative position of two intervals with integer
endpoints, we get a complex polynomial $P(z)$. The question of how
many roots are on the circle $\left|z\right|=1$ decide a number of
properties of the spectrum of an associated operator $P_{B}$. 
\begin{example}
\label{examples:ratinal vs irrational}Let $\Omega=[0,1]\cup[\alpha,\beta]$,
and fix the matrix $B=B(w,\phi,\psi,\theta)$. The spectrum $\Lambda_{B}$
of  $P_{B}$ is determined by 
\[
e\left(\theta-\phi+\left(\beta-\alpha\right)\lambda\right)=\frac{w\: e\left(\phi+\lambda\right)-1}{e\left(\phi+\lambda\right)-w}=e(g(t)).
\]
See Lemma \ref{Lemma:Magic-Equation} and \ref{Lemma:Mobius} for
details. That is, $\lambda\in\Lambda_{B}$ iff it is at the intersection
of the linear function $\theta-\phi+\left(\beta-\alpha\right)\lambda\;(\mbox{mod }1)$,
and the monotone decreasing function $g\left(t\right)\:(\mbox{mod 1})$.
Density of $\Lambda_{B}$ is reflected in the slope $(\beta-\alpha)$. 

Choose $w=\frac{1}{\sqrt{3}}$, $\theta=-1/4$, $\phi=-1/8$, $\alpha=3$,
and let $\beta=4,5,6,7$. Since $\beta-\alpha$ is rational, the spectrum
is periodic; see Figure \ref{fig:rat v.s. irrat 1}. On the other
hand, let $\beta=3+\{1,2,3,4\}\sqrt{2}$, so that $\beta-\alpha$
is irrational, and the spectrum is purely aperiodic; see Figure \ref{fig:rat v.s. irrat 2}.
\end{example}

\bigskip{}

\subsection{Examples}

What follows are examples \ref{Example:standard-sp} - \ref{example:(0,1)(3,6)};
and a summary of selected conclusions from these examples.

In Examples \ref{Example:standard-sp} and \ref{examples:(0,1)(2,3)},
we compute the boundary matrix $B$ for the case when $\alpha=2$,
$\beta=3$. The union of the corresponding two intervals offers a
spectral pair which is not a lattice tile.

Continuing with cases when the second interval has the form $[\alpha,\beta]$
with integer endpoints, we follow with a family of examples $\beta-\alpha=1$;
i.e., same length. They again are spectral pairs, but when $\alpha$
is odd, then they are also in the more restrictive class, the lattice
tiles. Indeed, for $\alpha$ odd, it turns out that there is then
both a lattice spectrum, and a non-lattice spectrum. See Example \ref{examples:(0,1)(3,4)}
which also shows that different choices of boundary matrix $B$, while
leading to the same spectrum, yet can result in different systems
of eigenfunctions. In Figures \ref{fig:(0,1)(2,3)} - \ref{fig:(0,1)(3,4)},
we illustrate determination of spectra from curve crossings. Examples
\ref{example:(0,1)(3,5)} - \ref{example:(0,1)(3,6)} illustrate computation
of spectrum($P_{B}$) for sets that are unions of pairs of intervals,
different length, and producing non-spectral pairs, i.e., not spectral
pairs for any choice of $\Lambda$.
\begin{example}
\label{Example:standard-sp} We know that $[0,1]\cup[2,3]$ is a spectral
set with spectrum $\mathbb{Z}\cup\left(\frac{1}{4}+\mathbb{Z}\right)$.
In this case $B=\left(\begin{array}{cc}
w & x\\
y & z
\end{array}\right)$ must satisfy the equations
\begin{align*}
\left(\begin{array}{cc}
w & x\\
y & z
\end{array}\right)\left(\begin{array}{c}
1\\
1
\end{array}\right) & =\left(\begin{array}{c}
1\\
1
\end{array}\right),\,\lambda\,\text{integer}\\
\left(\begin{array}{cc}
w & x\\
y & z
\end{array}\right)\left(\begin{array}{c}
i\\
-i
\end{array}\right) & =\left(\begin{array}{c}
1\\
-1
\end{array}\right),\,\lambda\,\text{integer plus 1/4}
\end{align*}
Hence 
\begin{equation}
B=\frac{1}{\sqrt{2}}\left(\begin{array}{cc}
1-i & 1+i\\
1+i & 1-i
\end{array}\right).\label{eq:B-temp}
\end{equation}

The boundary matrix $B$ in (\ref{eq:B-temp}) has parameters $w=\sqrt{2}/2$,
$\phi=-\frac{1}{8}$, $\psi=\frac{1}{8}$, and $\theta=-\frac{1}{4}$.
By (\ref{eq:PhaseIntersection (mix) 1}), $\lambda$ is in the spectrum
if and only if 
\[
e(\lambda-\frac{1}{8})=\frac{\sqrt{2}}{2}\cdot\frac{e(\lambda-\frac{1}{8})-\sqrt{2}}{e(\lambda-\frac{1}{8})-\frac{\sqrt{2}}{2}}.
\]
Since both sides of the above equation have modulus one, comparing
arguments shows that 
\[
\lambda-\frac{1}{8}=\pm\frac{1}{8}+\mathbb{Z}
\]
i.e., $\lambda=\{0,\frac{1}{4}\}+\mathbb{Z}=\mathbb{Z}\cup(\frac{1}{4}+\mathbb{Z})$
confirming the claim at the beginning. We revisit this example several
times below.
\end{example}

In the following examples we apply Corollary \ref{cor:Spectral-Set-Criterion}
in order to decide when $L^{2}(I_{1}\cup I_{2})$ has a Fourier basis,
i.e., when we get spectral pairs. In the examples below this amounts
to counting the number of root on the unit circle of a certain polynomial,
see also \cite{Ped96}.

\begin{figure}[H]
\begin{centering}
\includegraphics{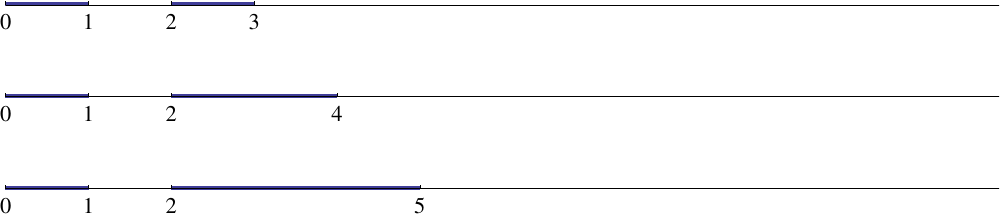}
\par\end{centering}

\caption{\label{fig:Class-A}Class $A$: $[0,1]\cup[2,2+k]$}
\end{figure}

\begin{figure}[H]
\begin{centering}
\includegraphics{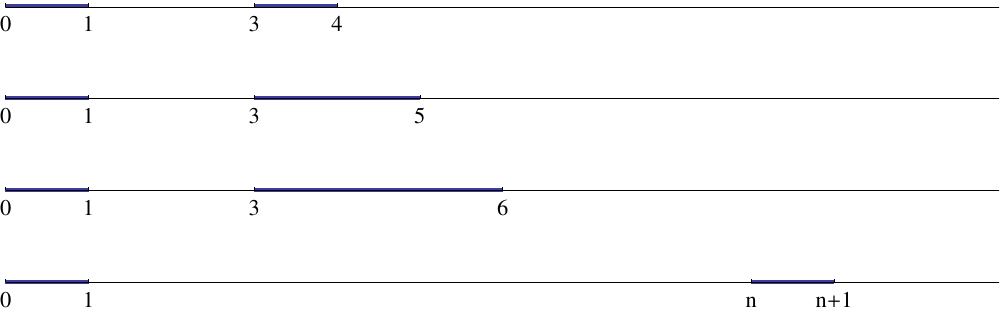}
\par\end{centering}

\caption{\label{fig:Class-B}Class $B$: $[0,1]\cup[3,3+k]$}
\end{figure}

We illustrate the results with two classes of examples, see Figures
\ref{fig:Class-A}, \ref{fig:Class-B}. Fix the parameters by $w=\frac{1}{\sqrt{2}}$,
$\theta=-1/4$, $\phi=-1/8$, $\psi=1/8$, and consider the two intervals
$I_{1},I_{2}$ given by $[0,1]\cup[2,2+k]$ in Class $A$, and $[0,1]\cup[3,3+k]$
in Class $B$, where $k=1,2,3\ldots$. Note that
\begin{enumerate}
\item The spectrum $\Lambda$ depends on $\beta-\alpha$ but not the location
of the interval $[\alpha,\beta]$;
\item The coefficient $a(\lambda)$ depends on all $6$ parameters. 
\end{enumerate}

\begin{example}
\label{examples:(0,1)(2,3)}Spectral pair: $[0,1]\cup[2,3]$ , $w=\frac{1}{\sqrt{2}}$,
$\theta=-1/4$, $\phi=-1/8$, $\psi=1/8$.\\
Spectrum 
\[
\Lambda=\{0,1/4\}+\mathbb{Z}=\mathbb{Z}\cup\left(1/4+\mathbb{Z}\right)
\]
Coefficients 
\[
a_{\lambda}=1,\forall\lambda\in\Lambda
\]
\begin{figure}[H]
\begin{centering}
\begin{tabular}{cc}
\includegraphics{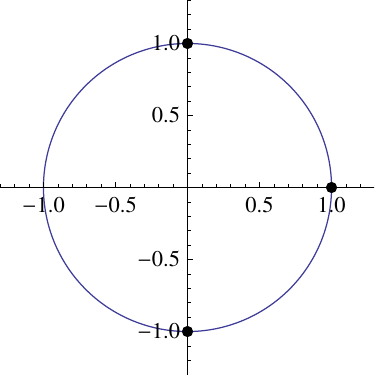} & \includegraphics{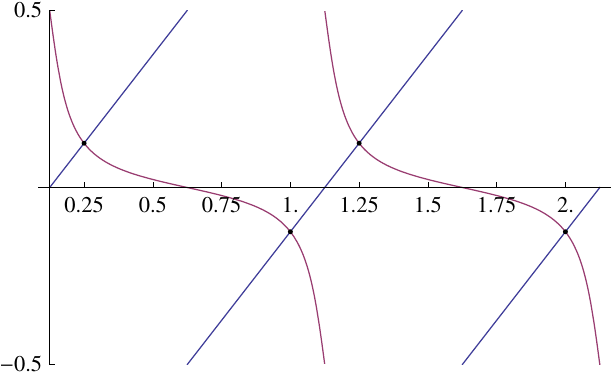}\tabularnewline
$p(z)=(z-1)(1+z^{2})$ & $\Lambda$\tabularnewline
\end{tabular}
\par\end{centering}

\caption{\label{fig:(0,1)(2,3)}}
\end{figure}

\end{example}

\begin{example}
\label{example:(0,1)(2,4)}Not spectral pair: $[0,1]\cup[2,4]$ ,
$w=\frac{1}{\sqrt{2}}$, $\theta=-1/4$, $\phi=-1/8$, $\psi=1/8$.\\
Spectrum 
\[
\Lambda\thickapprox\{0,0.1825,0.5675\}+\mathbb{Z}
\]
Coefficients 
\begin{eqnarray*}
a_{0} & = & 1\\
a_{0.1825} & \thickapprox & 1.61716-0.455719i\\
a_{0.5675} & \thickapprox & -0.367157+0.205719i
\end{eqnarray*}

\begin{figure}[H]
\begin{centering}
\begin{tabular}{cc}
\includegraphics{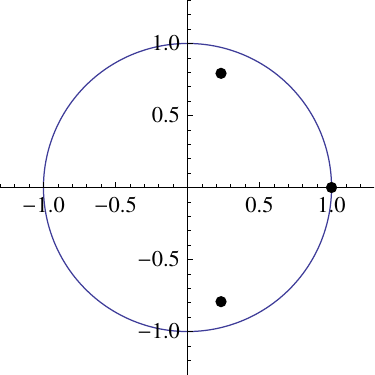} & \includegraphics{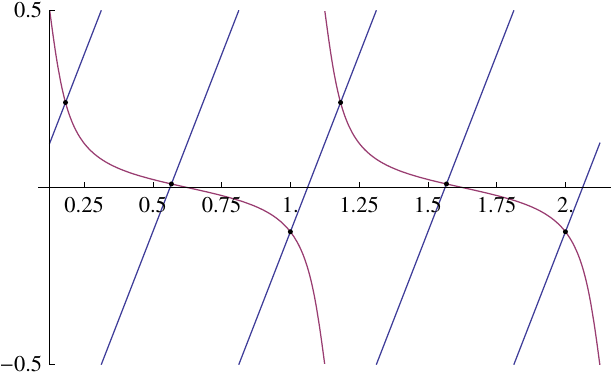}\tabularnewline
$p(z)=(z-1)(1+z^{2}(1+z))$ & $\Lambda$\tabularnewline
\end{tabular}
\par\end{centering}

\caption{\label{fig:(0,1)(2,4)}}
\end{figure}

\end{example}

\begin{example}
\label{example:(0,1)(2,5)}Not spectral pair: $[0,1]\cup[2,5]$ ,
$w=\frac{1}{\sqrt{2}}$, $\theta=-1/4$, $\phi=-1/8$, $\psi=1/8$.\\
Spectrum 
\[
\Lambda\thickapprox\{0,0.1550,0.3913,0.7037\}+\mathbb{Z}
\]
Coefficients
\begin{eqnarray*}
a_{0} & = & 1\\
a_{0.1550} & \thickapprox & 2.07306-0.464978i\\
a_{0.3913} & \thickapprox & 0.253715+0.489548i\\
a_{0.7037} & \thickapprox & -0.326779-0.27457i
\end{eqnarray*}

\begin{figure}[H]
\begin{centering}
\begin{tabular}{cc}
\includegraphics{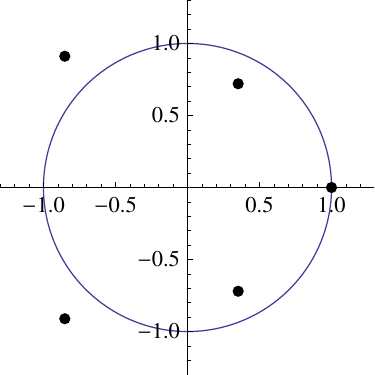} & \includegraphics{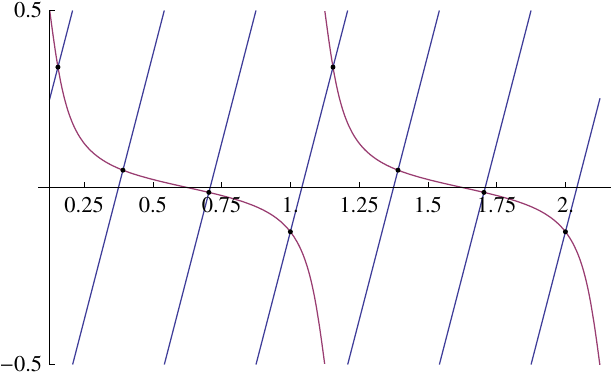}\tabularnewline
$p(z)=(z-1)(1+z^{2}(1+z+z^{2}))$ & $\Lambda$\tabularnewline
\end{tabular}
\par\end{centering}

\caption{\label{fig:(0,1)(2,5)}}
\end{figure}

\end{example}

\begin{example}
\label{examples:(0,1)(3,4)}Spectral pair: $[0,1]\cup[3,4]$ , $w=\frac{1}{\sqrt{2}}$,
$\theta=-1/4$, $\phi=-1/8$, $\psi=1/8$.\\
Spectrum 
\[
\Lambda=\{0,1/4\}+\mathbb{Z}=\mathbb{Z}\cup\left(1/4+\mathbb{Z}\right)
\]
But the coefficients are 
\[
a_{0}=1,\; a_{1/4}=i
\]
so the parameters above lead to non-spectral eigenfunctions. Illustrating
that different parameters $B=B(\theta,\psi,\phi,w)$ can lead to the
same spectrum yet to different eigenfunctions.

\begin{figure}[H]
\begin{centering}
\begin{tabular}{cc}
\includegraphics{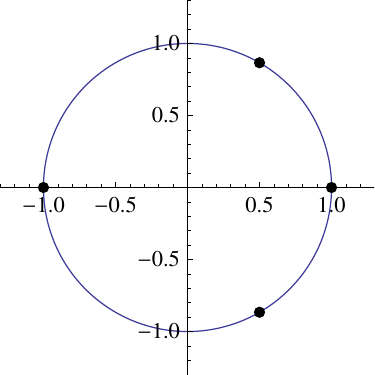} & \includegraphics{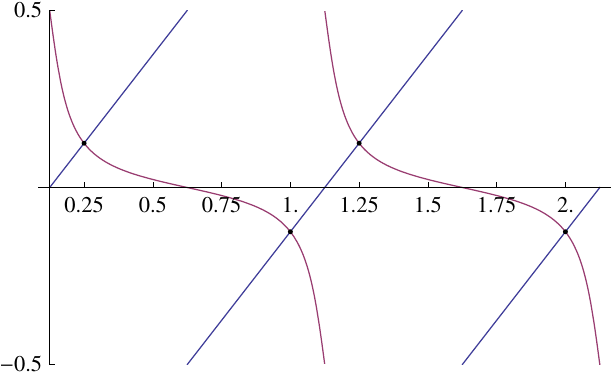}\tabularnewline
$p(z)=(z-1)(1+z^{3})$ & $\Lambda$\tabularnewline
\end{tabular}
\par\end{centering}

\caption{\label{fig:(0,1)(3,4)}}
\end{figure}

\end{example}

\begin{example}
\label{example:(0,1)(3,5)}Not spectral pair: $[0,1]\cup[3,5]$ ,
$w=\frac{1}{\sqrt{2}}$, $\theta=-1/4$, $\phi=-1/8$, $\psi=1/8$.
\\
Spectrum 
\[
\Lambda\thickapprox\{0,0.1825,0.5675\}+\mathbb{Z}
\]
Coefficients
\begin{eqnarray*}
a_{0} & = & 1\\
a_{0.1825} & \thickapprox & 1.08072+1.28644i\\
a_{0.5675} & \thickapprox & 0.419281-0.0364378i
\end{eqnarray*}

\begin{figure}[H]
\begin{centering}
\begin{tabular}{cc}
\includegraphics{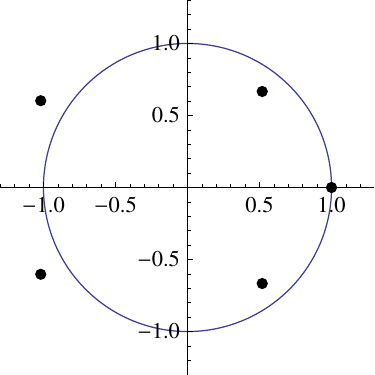} & \includegraphics{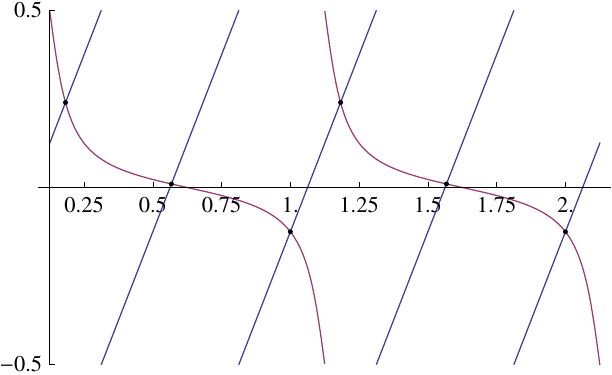}\tabularnewline
$p(z)=(z-1)(1+z^{3}(1+z))$ & $\Lambda$\tabularnewline
\end{tabular}
\par\end{centering}

\caption{\label{fig:(0,1)(3,5)}}
\end{figure}

\end{example}

\begin{example}
\label{example:(0,1)(3,6)}Not spectral pair: $[0,1]\cup[3,6]$ ,
$w=\frac{1}{\sqrt{2}}$, $\theta=-1/4$, $\phi=-1/8$, $\psi=1/8$,
$p(z)=(z-1)(1+z^{3}(1+z+z^{2}))$.\\
Spectrum 
\[
\Lambda\thickapprox\{0,0.1550,0.3913,0.7037\}+\mathbb{Z}
\]
Coefficients
\begin{eqnarray*}
a_{0} & = & 1\\
a_{0.1550} & \thickapprox & 1.55016+1.45286i\\
a_{0.3913} & \thickapprox & -0.505763-0.219617i\\
a_{0.7037} & \thickapprox & -0.1694+0.391761i
\end{eqnarray*}

\begin{figure}[H]
\begin{centering}
\begin{tabular}{cc}
\includegraphics{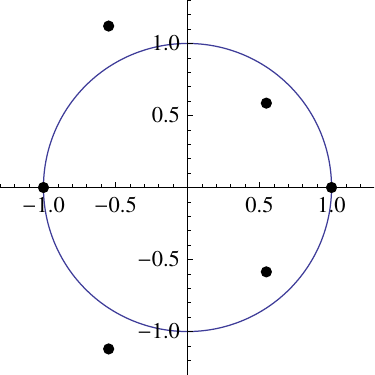} & \includegraphics{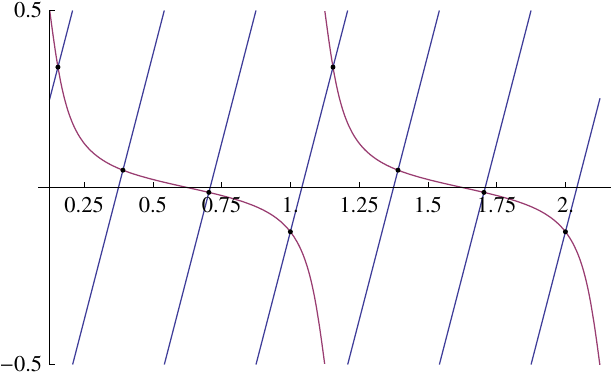}\tabularnewline
$p(z)=(z-1)(1+z^{3}(1+z))$ & $\Lambda$\tabularnewline
\end{tabular}
\par\end{centering}

\caption{\label{fig:(0,1)(3,6)}}
\end{figure}

\end{example}

\section{Questions}

Motivated by connections to problems from mathematical physics, we
include a brief list of questions arising naturally in connection
with properties of spectrum($P_{B}$) for the selfadjoint operators
$P_{B}$ considered here. Even though \cite{Rob71} deals with second
order operators, there are a number of parallels, for example both
in \cite{Rob71} and here, we are dealing with deficiency indices
$(2,2)$. And in both cases, we get that, when the boundary points
are fixed, then the individual eigenvalues $\lambda_{n}(B)$ from
the list spectrum($P_{B}$), for each $n$, depend continuously on
$B$. 
\begin{enumerate}
\item Robinson \cite{Rob71}: one interval, but a second order operator,
indices $(2,2),$ and a formula where the spectrum is a sequence of
points in $\lambda$ (on the real line) arising by graph intersections.
Our case: two intervals, and a first order operator, and the same
conclusion as in Robinson, only with graphs of a different pair of
functions. Open question: What can be said about \textquotedblleft{}similar\textquotedblright{}
of extensions of Hermitian operators with indices $(2,2)?$ 
\item Which self-adjoint restrictions of $d^{2}/dx^{2}$ are not $P^{2}$
for some self-adjoint restriction $P$ of $id/dx$? 
\item Obtain a formula for the dynamics, i.e., description of $t\to e^{itP_{B}}f(x)?$ 
\item Corollary \ref{cor:Spectral-separation} and plots of $\lambda_{n}$
as a function of $n,$ suggests it might be interesting to obtain
a quantitative estimate, depending on the other parameters, on $\delta=\delta_{w}?$ 
\item More than two intervals. Some results in the paper do carry over to
sets formed as union of $n$ intervals, but not our detailed classifications.
With $n$ intervals, the preliminaries in Sections \ref{sec:Introduction}
- \ref{sec:Momentum-Operators} still hold: the deficiency indices
will be $(n,n)$. The boundary form in Lemma \ref{lem:vN def-space}
will also have a matrix representation. Corollary \ref{cor:matrix-realization}
will still hold; but instead there will be $n$ left interval endpoints,
and a corresponding set of right-hand side interval endpoints. As
a result, the boundary matrix $B$ will now be an $n$ by $n$ complex
unitary matrix; every $B$ in the group $U(n)$ will define a selfadjoint
operator $P_{B}$ . The group $U(n)$ allows an Iwasawa decomposition,
and we will get an analogue of $w$ in (\ref{eq:2-by-2 Unitary}),
only now it will be an positive selfadjoint matrix $W$, satisfying
$0\leq W\leq I$, referring to the usual order of selfadjoint matrices.
One will also have a matrix version of the M\"{o}bius transform in
(\ref{eq:Mobius}) in Lemma \ref{Lemma:Mobius}. There is a theory
of matrix operations, making use of fractional linear matrix operations.
There is even a matrix version of the Poisson kernel (\ref{eq:w-possion}).
But our detailed classification in Theorem \ref{thm:0<w<1} will not
carry over. Nonetheless, each selfadjoint operator $P_{B}$ will still
have a pure point spectrum, but now with a much more subtle multiplicity
configuration. And, for the general multiple-interval case, it would
be difficult to get a classification of all the spectral pairs anything
close to what we have in our theorems in Section \ref{sec:Spectral-Pairs}
above.
\end{enumerate}

\section*{Acknowledgments}

The co-authors, some or all, had helpful conversations with many colleagues,
and wish to thank especially Professors Daniel Alpay, Ilwoo Cho, Dorin
Dutkay, Alex Iosevich, Paul Muhly, and Yang Wang. And going back in
time, Bent Fuglede (PJ, SP), and Derek Robinson (PJ).

\bibliographystyle{alpha}
\bibliography{Intervals}

\begin{thebibliography}{BAAS11}

\bibitem[BAAS11]{BAAS11}
S.~Behnia, A.~Akhavan, A.~Akhshani, and A.~Samsudin.
\newblock A novel dynamic model of pseudo random number generator.
\newblock {\em J. Comput. Appl. Math.}, 235(12):3455--3463, 2011.

\bibitem[Bar49]{Ba49}
V.~Bargmann.
\newblock On the connection between phase shifts and scattering potential.
\newblock {\em Rev. Modern Physics}, 21:488--493, 1949.

\bibitem[BH08]{BH08}
Horst Behncke and D.~B. Hinton.
\newblock Eigenfunctions, deficiency indices and spectra of odd-order
  differential operators.
\newblock {\em Proc. Lond. Math. Soc. (3)}, 97(2):425--449, 2008.

\bibitem[BV05]{BV05}
Pallav~Kumar Baruah and M.~Venkatesulu.
\newblock Deficiency indices of a differential operator satisfying certain
  matching interface conditions.
\newblock {\em Electron. J. Differential Equations}, pages No. 38, 9 pp.
  (electronic), 2005.

\bibitem[DJ07]{DJ07}
Dorin~Ervin Dutkay and Palle E.~T. Jorgensen.
\newblock Fourier frequencies in affine iterated function systems.
\newblock {\em J. Funct. Anal.}, 247(1):110--137, 2007.

\bibitem[dO09]{dO09}
C{\'e}sar~R. de~Oliveira.
\newblock {\em Intermediate spectral theory and quantum dynamics}, volume~54 of
  {\em Progress in Mathematical Physics}.
\newblock Birkh\"auser Verlag, Basel, 2009.

\bibitem[DS88]{DS88b}
Nelson Dunford and Jacob~T. Schwartz.
\newblock {\em Linear operators. {P}art {II}}.
\newblock Wiley Classics Library. John Wiley \& Sons Inc., New York, 1988.
\newblock Spectral theory. Selfadjoint operators in Hilbert space, With the
  assistance of William G. Bade and Robert G. Bartle, Reprint of the 1963
  original, A Wiley-Interscience Publication.

\bibitem[Fug74]{Fu74}
Bent Fuglede.
\newblock Commuting self-adjoint partial differential operators and a group
  theoretic problem.
\newblock {\em J. Functional Analysis}, 16:101--121, 1974.

\bibitem[Gil72]{Gil72}
Richard~C. Gilbert.
\newblock Spectral representation of selfadjoint extensions of a symmetric
  operator.
\newblock {\em Rocky Mountain J. Math.}, 2(1):75--96, 1972.

\bibitem[IKT03]{IKT03}
Alex Iosevich, Nets Katz, and Terence Tao.
\newblock The {F}uglede spectral conjecture holds for convex planar domains.
\newblock {\em Math. Res. Lett.}, 10(5-6):559--569, 2003.

\bibitem[IP98]{IP98}
Alex Iosevich and Steen Pedersen.
\newblock Spectral and tiling properties of the unit cube.
\newblock {\em Internat. Math. Res. Notices}, (16):819--828, 1998.

\bibitem[J{\o}r81]{Jo81}
Palle E.~T. J{\o}rgensen.
\newblock A uniqueness theorem for the {H}eisenberg-{W}eyl commutation
  relations with nonselfadjoint position operator.
\newblock {\em Amer. J. Math.}, 103(2):273--287, 1981.

\bibitem[J{\o}r82]{Jo82}
Palle E.~T. J{\o}rgensen.
\newblock Spectral theory of finite volume domains in {${\bf R}^{n}$}.
\newblock {\em Adv. in Math.}, 44(2):105--120, 1982.

\bibitem[JP98]{JP98}
Palle E.~T. Jorgensen and Steen Pedersen.
\newblock Dense analytic subspaces in fractal {$L^2$}-spaces.
\newblock {\em J. Anal. Math.}, 75:185--228, 1998.

\bibitem[JP99]{JP99}
Palle E.~T. Jorgensen and Steen Pedersen.
\newblock Spectral pairs in {C}artesian coordinates.
\newblock {\em J. Fourier Anal. Appl.}, 5(4):285--302, 1999.

\bibitem[JP00]{JP00}
Palle E.~T. Jorgensen and Steen Pedersen.
\newblock Commuting self-adjoint extensions of symmetric operators defined from
  the partial derivatives.
\newblock {\em J. Math. Phys.}, 41(12):8263--8278, 2000.

\bibitem[KM06]{KM06}
Mihail~N. Kolountzakis and M{\'a}t{\'e} Matolcsi.
\newblock Complex {H}adamard matrices and the spectral set conjecture.
\newblock {\em Collect. Math.}, (Vol. Extra):281--291, 2006.

\bibitem[Kre49]{Kr49}
M.~G. Kre{\u\i}n.
\newblock The fundamental propositions of the theory of representations of
  {H}ermitian operators with deficiency index {$(m,m)$}.
\newblock {\em Ukrain. Mat. \v Zurnal}, 1(2):3--66, 1949.

\bibitem[{\L}ab01]{Laba01}
I.~{\L}aba.
\newblock Fuglede's conjecture for a union of two intervals.
\newblock {\em Proc. Amer. Math. Soc.}, 129(10):2965--2972 (electronic), 2001.

\bibitem[{\L}ab02]{La02}
I.~{\L}aba.
\newblock The spectral set conjecture and multiplicative properties of roots of
  polynomials.
\newblock {\em J. London Math. Soc. (2)}, 65(3):661--671, 2002.

\bibitem[Lan67]{Landau}
H.~J. Landau.
\newblock Necessary density conditions for sampling and interpolation of
  certain entire functions.
\newblock {\em Acta Math.}, 117:37--52, 1967.

\bibitem[LP68]{LP68}
P.~D. Lax and R.~S. Phillips.
\newblock Scattering theory.
\newblock In {\em Proc. {I}nternat. {C}ongr. {M}ath. ({M}oscow, 1966)}, pages
  542--545. Izdat. ``Mir'', Moscow, 1968.

\bibitem[LW97]{LW97}
Jeffrey~C. Lagarias and Yang Wang.
\newblock Spectral sets and factorizations of finite abelian groups.
\newblock {\em J. Funct. Anal.}, 145(1):73--98, 1997.

\bibitem[Mik04]{Mik04}
V.~A. Mikha{\u\i}lets.
\newblock The general spectrum of a family of selfadjoint extensions.
\newblock {\em Dopov. Nats. Akad. Nauk Ukr. Mat. Prirodozn. Tekh. Nauki},
  (1):18--21, 2004.

\bibitem[Min04]{Min04}
V.~S. Mineev.
\newblock Physics of selfadjoint extensions: the one-dimensional scattering
  problem for {C}oulomb potential.
\newblock {\em Teoret. Mat. Fiz.}, 140(2):310--328, 2004.

\bibitem[Naz08]{Naz08}
S.~A. Nazarov.
\newblock Selfadjoint extensions of the operator of the {D}irichlet problem in
  a three-dimensional domain with an edge.
\newblock {\em Sib. Zh. Ind. Mat.}, 11(1):80--95, 2008.

\bibitem[Oro05]{Oro05}
Yu.~B. Orochko.
\newblock Deficiency indices of an even-order one-term symmetric differential
  operator that degenerates inside an interval.
\newblock {\em Mat. Sb.}, 196(5):53--82, 2005.

\bibitem[Ped96]{Ped96}
Steen Pedersen.
\newblock Spectral sets whose spectrum is a lattice with a base.
\newblock {\em J. Funct. Anal.}, 141(2):496--509, 1996.

\bibitem[Ped04]{Pe04}
Steen Pedersen.
\newblock The dual spectral set conjecture.
\newblock {\em Proc. Amer. Math. Soc.}, 132:2095--2101, 2004.

\bibitem[PR76]{PoRa76}
Robert~T. Powers and Charles Radin.
\newblock Average boundary conditions in {C}auchy problems.
\newblock {\em J. Functional Analysis}, 23(1):23--32, 1976.

\bibitem[PW01]{PW01}
Steen Pedersen and Yang Wang.
\newblock Universal spectra, universal tiling sets and the spectral set
  conjecture.
\newblock {\em Math. Scand.}, 88(2):246--256, 2001.

\bibitem[Rob71]{Rob71}
Derek~W. Robinson.
\newblock {\em The thermodynamic pressure in quantum statistical mechanics}.
\newblock Springer-Verlag, Berlin, 1971.
\newblock Lecture Notes in Physics, Vol. 9.

\bibitem[Sad06]{Sad06}
I.~V. Sadovnichaya.
\newblock A new estimate for the spectral function of a selfadjoint extension
  in {$L^2(\Bbb R)$} of the {S}turm-{L}iouville operator with a uniformly
  locally integrable potential.
\newblock {\em Differ. Uravn.}, 42(2):188--201, 286, 2006.

\bibitem[Sak97]{Sak97}
L.~A. Sakhnovich.
\newblock Deficiency indices of a system of first-order differential equations.
\newblock {\em Sibirsk. Mat. Zh.}, 38(6):1360--1361, iii, 1997.

\bibitem[{\v{S}}mu74]{Smu74}
Ju.~L. {\v{S}}mul{\cprime}jan.
\newblock Closed {H}ermitian operators and their selfadjoint extensions.
\newblock {\em Mat. Sb. (N.S.)}, 93(135):155--169, 325, 1974.

\bibitem[ST10]{ST10}
Luis~O. Silva and Julio~H. Toloza.
\newblock On the spectral characterization of entire operators with deficiency
  indices {$(1,1)$}.
\newblock {\em J. Math. Anal. Appl.}, 367(2):360--373, 2010.

\bibitem[Tao04]{Tao04}
Terence Tao.
\newblock Fuglede's conjecture is false in 5 and higher dimensions.
\newblock {\em Math. Res. Lett.}, 11(2-3):251--258, 2004.

\bibitem[Vas07]{Vas07}
F.-H. Vasilescu.
\newblock Existence of the smallest selfadjoint extension.
\newblock In {\em Perspectives in operator theory}, volume~75 of {\em Banach
  Center Publ.}, pages 323--326. Polish Acad. Sci., Warsaw, 2007.

\bibitem[VGT08]{VGT08}
B.~L. Voronov, D.~M. Gitman, and I.~V. Tyutin.
\newblock Construction of quantum observables and the theory of selfadjoint
  extensions of symmetric operators. {III}. {S}elfadjoint boundary conditions.
\newblock {\em Izv. Vyssh. Uchebn. Zaved. Fiz.}, 51(2):3--43, 2008.

\bibitem[vN32]{vNeu32}
J.~von Neumann.
\newblock \"{U}ber adjungierte {F}unktionaloperatoren.
\newblock {\em Ann. of Math. (2)}, 33(2):294--310, 1932.

\end{thebibliography}

\end{document}